\documentclass{agtart_a}
\pdfoutput=1

\usepackage{graphicx}
\usepackage{pinlabel}

\title{Holomorphic discs and sutured manifolds}

\author{Andr\'as Juh\'asz}
\givenname{Andr\'as}
\surname{Juh\'asz}
\address{Department of Mathematics\\
Princeton University\\\newline
Princeton, NJ 08544\\USA}
\email{ajuhasz@math.princeton.edu}
\urladdr{}

\volumenumber{6}
\issuenumber{}
\publicationyear{2006}
\papernumber{52}
\startpage{1429}
\endpage{1457}

\doi{}
\MR{}
\Zbl{}


\received{18 March 2006}
\revised{}
\accepted{23 July 2006}
\published{4 October 2006}
\publishedonline{4 October 2006}
\proposed{}
\seconded{}
\corresponding{}
\editor{}
\version{}

\arxivreference{math.GT/0601443}

\makeop{Tw}
\makeop{SFH}
\makeop{HF}
\makeop{HFL}
\makeop{HFK}
\def\tsty{\textstyle}

\makeatletter
\def\cnewtheorem#1[#2]#3{\newtheorem{#1}{#3}[section]
\expandafter\let\csname c@#1\endcsname\c@thm}
\makeatother

\newtheorem{thm}{Theorem}[section]
\cnewtheorem{cor}[thm]{Corollary}
\cnewtheorem{lem}[thm]{Lemma}
\cnewtheorem{prop}[thm]{Proposition}
\cnewtheorem{conj}[thm]{Conjecture}
\theoremstyle{definition}
\cnewtheorem{defn}[thm]{Definition}
\makeautorefname{defn}{Definition}
\cnewtheorem{probl}[thm]{Problem}
\cnewtheorem{qn}[thm]{Question}
\theoremstyle{remark}
\cnewtheorem{rem}[thm]{Remark}
\cnewtheorem{exm}[thm]{Example}
\makeautorefname{exm}{Example}
\numberwithin{equation}{section}
\cnewtheorem{note}[thm]{Notation}

\begin{document}

\begin{asciiabstract}
In this paper we construct a Floer-homology invariant for a natural
and wide class of sutured manifolds that we call balanced. This
generalizes the Heegaard Floer hat theory of closed three-manifolds
and links. Our invariant is unchanged under product decompositions
and is zero for nontaut sutured manifolds. As an application, an
invariant of Seifert surfaces is given and is computed in a few
interesting cases.
\end{asciiabstract}

\begin{abstract}
In this paper we construct a Floer-homology invariant for a natural
and wide class of sutured manifolds that we call balanced. This
generalizes the Heegaard Floer hat theory of closed three-manifolds
and links. Our invariant is unchanged under product decompositions
and is zero for nontaut sutured manifolds. As an application, an
invariant of Seifert surfaces is given and is computed in a few
interesting cases.
\end{abstract}

\maketitle
\section{Introduction}

In Ozsv{\'a}th and Szab{\'o} \cite{OSz} a Floer homology invariant was defined for closed
oriented $3$--manifolds. This theory was extended to knots by
Ozsv{\'a}th and Szab{\'o} \cite{OSz3} and Rasmussen \cite{Ras} and recently to links 
again by Ozsv{\'a}th and Szab{\'o} \cite{OSz2}.
Motivated by a conjecture that knot Floer homology detects fibred
knots (\fullref{conj:2}, originally proposed in \cite{OSz4})
and a characterization of fibred knots by Gabai \cite{Gabai2}, we
extend Heegaard Floer hat theory to a class of sutured manifolds
that we call balanced (\fullref{defn:2}). This theory
provides us with a new invariant that we call sutured Floer
homology, in short, $\SFH$. In particular, for every closed oriented
$3$--manifold $Y$ and every link $L \subset Y$ we construct balanced
sutured manifolds $Y(1)$ and $Y(L)$ such that $\widehat{\HF}(Y) =
\SFH(Y(1))$ and $\widehat{\HFL}(L) = \SFH(Y(L)) \otimes \mathbb{Z}_2$.
Any group $\SFH(M, \gamma)$ decomposes into a direct sum along
relative $\text{Spin}^c$ structures on the sutured manifold
$(M,\gamma)$ and each summand possesses a relative grading.

To construct the invariant we define the notion of a balanced
Heegaard diagram (\fullref{defn:5}), which consists of a
compact surface $\Sigma$ with no closed components and sets of
curves $\boldsymbol{\alpha}$ and $\boldsymbol{\beta}$ of the same
cardinality $d$ that are also linearly independent in $H_1(\Sigma;
\mathbb{Q})$. These data provide the input for the usual
construction of Lagrangian Floer homology applied to
$\mathbb{T}_{\alpha}, \mathbb{T}_{\beta} \subset
\text{Sym}^d(\Sigma)$.

The invariant that we have constructed is unchanged under product
decompositions of sutured manifolds (\fullref{lem:9}) and is zero
for nontaut sutured manifolds (\fullref{prop:8}). In the
last chapter we assign to every Seifert surface $R \subset S^3$ a
sutured manifold $S^3(R)$ and we compute $\SFH(S^3(R))$ in a few
cases. These computations indicate a relationship between the top
nonzero term of knot Floer homology and sutured Floer homology of
the sutured manifold obtained from a minimal genus Seifert surface.
This relationship is the subject of \fullref{conj:1}.

\paragraph{Acknowledgements}
I would like to thank Professor Zolt\'an Szab\'o for leading me to
the idea of sutured Floer homology and for his support during the
course of this work. I would also like to thank Yi Ni for the
helpful discussions about the topic, for thoroughly reading the
first version of this paper and for proving \fullref{prop:8}. I am grateful to the referee for carefully reading the
manuscript and making several useful remarks. This research was
partially supported by OTKA grant no. T49449.

\section{Heegaard diagrams of sutured manifolds}

First we recall the notion of a sutured manifold as defined by Gabai
\cite{Gabai}.

\begin{defn} \label{defn:1}
A \emph{sutured manifold\/} $(M,\gamma)$ is a compact oriented
$3$--manifold $M$ with boundary together with a set $\gamma \subset
\partial M$ of pairwise disjoint annuli $A(\gamma)$ and tori
$T(\gamma)$. Furthermore, the interior of each component of
$A(\gamma)$ contains a \emph{suture\/}, ie, a homologically
nontrivial oriented simple closed curve. We denote the union of the
sutures by $s(\gamma)$.

Finally every component of $R(\gamma)=\partial M \setminus
\text{Int}(\gamma)$ is oriented. Define $R_+(\gamma)$ (or
$R_-(\gamma)$) to be those components of $\partial M \setminus
\text{Int}(\gamma)$ whose normal vectors point out of (into) $M$.
The orientation on $R(\gamma)$ must be coherent with respect to
$s(\gamma)$, ie, if $\delta$ is a component of $\partial
R(\gamma)$ and is given the boundary orientation, then $\delta$ must
represent the same homology class in $H_1(\gamma)$ as some suture.
\end{defn}

In this paper we will restrict our attention to a special class of
sutured manifolds.

\begin{defn} \label{defn:2}
A \emph{balanced sutured manifold\/} is a sutured manifold
$(M,\gamma)$ such that M has no closed components,
$\chi(R_+(\gamma))=\chi(R_-(\gamma))$, and the map from $\pi_0(A(\gamma))$
to $\pi_0(\partial M)$ is surjective.
\end{defn}

Note that the last condition implies that for a balanced sutured
manifold $T(\gamma)=\emptyset$. A balanced sutured manifold is
completely determined by $M$ and $s(\gamma)$. Therefore, one can
view $\gamma$ as a set of thick oriented curves in $\partial M$
where such curves induce the orientations on $\partial M \setminus
\text{Int}(\gamma)$. Now we list a few important examples of
balanced sutured manifolds.

\begin{exm} \label{ex:1}
Let $Y$ be a closed connected oriented $3$--manifold and we are also
given pairwise disjoint closed $3$--balls $B_1, \dots, B_k \subset Y$. 
For $1 \le i \le k$ choose an oriented simple closed curve $s_i
\subset
\partial B_i$ together with a regular neighborhood $\gamma_i = N(s_i)$.
If $$M = Y \setminus \tsty\bigcup_{i=1}^k \text{Int}(B_i)\quad\text{and}\quad\gamma =
\tsty\bigcup_{i=1}^k \gamma_i$$ then the pair $(M, \gamma)$ defines a
balanced sutured manifold with sutures $s(\gamma)=\bigcup_{i=1}^k
s_i$. The sutured manifold $(M,\gamma)$ only depends on $Y$ and $k$, 
we denote it by $Y(k)$. Note that $Y(k)$ uniquely determines $Y$. 

If $(N,\nu)$ is a connected balanced sutured manifold then let
$N(k)$ denote the connected sum $(N, \nu)\#S^3(k)$. This is also a
balanced sutured manifold.
\end{exm}

\begin{exm} \label{ex:2}
Let $L \subset Y$ be a link of $k$ components in a closed connected
oriented $3$--manifold $Y$. Choose a closed regular neighborhood $N(L)$
of $L$. For every component $L_i$ of $L$ ($1 \le i \le k$) take two
meridians $s_i$ and $s_{\smash{i}}'$ of $L_i$ oppositely oriented, that is,
$[s_i] = -[s_{\smash{i}}']$ in $H_1(\partial N(L_i);\mathbb{Z})$. Choose
regular neighborhoods $\gamma_i = N(s_i)$ and $\gamma_{\smash{i}}' = N(s_{\smash{i}}')$
in $\partial N(L_i)$ and let $\gamma = \bigcup_{i=1}^{\smash{k}} (\gamma_i
\cup \gamma_{\smash{i}}')$; furthermore let $M = Y \setminus \bigcup_{i=1}^{\smash{k}}
\text{Int}(N(L_i))$. This way we obtain a balanced sutured manifold
$(M,\gamma)$. We can reconstruct $L$ from $(M,\gamma)$ using Dehn
filling as follows. For each component $T^2_i$ of $\partial M$ glue
in a solid torus $S^1 \times D^2$ so that $\{1\} \times \partial
D^2$ maps to one component of $s(\gamma) \cap T^2_i$, let $L_i$ be
the image of $S^1 \times \{0\}$. Note that if we choose the other
component of $s(\gamma) \cap T^2_i$ only the orientation of $L_i$
changes, and choosing different images for the longitude $S^1 \times
\{1\}$ corresponds to choosing different framings of $L_i$. 

$(M,\gamma)$ is uniquely determined by the link $L;$ let us use the
notation $Y(L)$ for the sutured manifold $(M,\gamma)$. We saw above
that $Y(L)$ uniquely determines $L$. If in addition we fix an
ordering of the components of $s(\gamma) \cap T^2_i$  (ie, we
distinguish between $s_i$ and $s_i'$) we uniquely define an
orientation of $L$. 
\end{exm}

The following two examples can be found in \cite{Gabai2}.

\begin{exm} \label{ex:3}
Let $R$ be a compact oriented surface with no closed components.
Then there is an induced orientation on $\partial R$. Let $M=R
\times I$, define $\gamma =\partial R \times I$, finally put
$s(\gamma) =
\partial R \times \{1/2\}$. The balanced sutured manifold
$(M,\gamma)$ obtained by this construction is called a \emph{product
sutured manifold}.
\end{exm}

\begin{exm} \label{ex:4}
Let $Y$ be a closed connected oriented $3$--manifold and let $R \subset
Y$ be a compact oriented surface with no closed components. We
define $Y(R)=(M,\gamma)$ to be the sutured manifold where $M = Y
\setminus \text{Int}(R \times I)$, the suture $\gamma =
\partial R \times I$ and $s(\gamma) = \partial R \times
\{1/2\}$. Then $Y(R)$ is balanced.
\end{exm}

Next we introduce sutured Heegaard diagrams. They generalize
Heegaard diagrams of closed $3$--manifolds so that we can also describe
sutured manifolds.

\begin{defn} \label{defn:3}
A \emph{sutured Heegaard diagram\/} is a tuple $( \Sigma,
\boldsymbol{\alpha}, \boldsymbol{\beta})$, where $\Sigma$ is a
compact oriented surface with boundary and $\boldsymbol{\alpha} =
\{\,\alpha_1, \dots, \alpha_m\,\}$ and $\boldsymbol{\beta} =
\{\,\beta_1, \dots, \beta_n\,\}$ are two sets of pairwise disjoint
simple closed curves in $\text{Int}(\Sigma)$. 
\end{defn}

\begin{defn} \label{defn:4}
Every sutured Heegaard diagram $( \Sigma, \boldsymbol{\alpha},
\boldsymbol{\beta})$ uniquely \emph{defines\/} a sutured manifold $(M,
\gamma)$ using the following construction.

Let $M$ be the $3$--manifold obtained from $\Sigma \times I$ by
attaching $3$--dimensional $2$--handles along the curves $\alpha_i \times
\{0\}$ and $\beta_j \times \{1\}$ for $i=1, \dots, m$ and $j=1,
\dots, n$. The sutures are defined by taking $\gamma = \partial M
\times I$ and $s(\gamma)=
\partial M \times \{1/2\}$. 
\end{defn}

\begin{prop} \label{prop:1}
If $(M,\gamma)$ is defined by $(\Sigma, \boldsymbol{\alpha},
\boldsymbol{\beta})$ then $(M,\gamma)$ is balanced if and only if
$|\boldsymbol{\alpha}| = |\boldsymbol{\beta}|$ and the maps
$\pi_0(\partial \Sigma) \to \pi_0(\Sigma \setminus \bigcup
\boldsymbol{\alpha})$ and $\pi_0(\partial \Sigma) \to \pi_0(\Sigma
\setminus \bigcup \boldsymbol{\beta})$ are surjective. The second
condition is equivalent to saying that $\Sigma$ has no closed
components and the elements of $\boldsymbol{\alpha}$ and
$\boldsymbol{\beta}$ are both linearly independent in $H_1(\Sigma;
\mathbb{Q})$. 
\end{prop}

\begin{proof}
Since adding a $2$--handle increases the Euler characteristics of the
boundary by 2 (the boundary undergoes surgery along the attaching
circle) we get the equalities $\chi(R_+(\gamma)) = \chi(\Sigma)+2m$ and
$\chi(R_-(\gamma)) = \chi(\Sigma)+2n$. Thus 
$\chi(R_+(\gamma)) = \chi(R_-(\gamma))$ if and only if
$|\boldsymbol{\alpha}| = |\boldsymbol{\beta}|$. 

Note that every component of $\partial M$ contains a suture exactly
when $R_{-}(\gamma)$ and $R_+(\gamma)$ have no closed components.
Since $R_-(\gamma)$ is obtained from $\Sigma$ by performing surgery
along $\boldsymbol{\alpha}$, components of $\Sigma \setminus \bigcup
\boldsymbol{\alpha}$ naturally correspond to components of
$R_-(\gamma)$. Thus a component of $\Sigma \setminus \bigcup
\boldsymbol{\alpha}$ contains a component of $\partial \Sigma$ if
and only if the corresponding component of $R_-(\gamma)$ has
nonempty boundary. So $R_-(\gamma)$ has no closed components if and
only if the map $\pi_0(\partial \Sigma) \to \pi_0(\Sigma \setminus
\bigcup \boldsymbol{\alpha})$ is surjective. A similar argument can
be used for $R_+(\gamma)$. 

The last statement follows from \fullref{lem:1}.
\end{proof}

\begin{lem} \label{lem:1}
Let $\Sigma$ be a compact oriented surface with boundary and let
$\alpha \subset \text{Int}(\Sigma)$ be a one-dimensional submanifold
of $\Sigma$. Then the map $\pi_0(\partial \Sigma) \to \pi_0(\Sigma
\setminus \alpha)$ is injective if and only if $\Sigma$ has no
closed components and the components of $\alpha$ are linearly
independent in $H_1(\Sigma; \mathbb{Q})$. 
\end{lem}

\begin{proof}
In this proof every homology group is to be considered with
coefficients in $\mathbb{Q}$. The components of $\alpha$ are
linearly independent in $H_1(\Sigma)$ exactly when the map $i_*
\co H_1(\alpha) \to H_1(\Sigma)$ induced by the embedding $i
\co \alpha \hookrightarrow \Sigma$ is injective. Look at the
following portion of the long exact sequence of the pair $(\Sigma,
\alpha):$
$$0 \to H_2(\Sigma) \to H_2(\Sigma, \alpha) \to H_1(\alpha)
\xrightarrow{i_*} H_1(\Sigma).$$ Then we see that $H_2(\Sigma,
\alpha) \approx H_2(\Sigma) \oplus \ker(i_*)$. Note that
$H_2(\Sigma) = 0$ precisely when $\Sigma$ has no closed components.
Let $N(\alpha)$ be a closed regular neighborhood of $\alpha$. Then
by excision
$$H_2(\Sigma, \alpha) \approx H_2(\Sigma, N(\alpha)) \approx
H_2(\Sigma \setminus \text{Int}(N(\alpha)), \partial N(\alpha))
\approx \bigoplus_C H_2(C, \partial N(\alpha) \cap C),$$ where $C$
runs over the components of $\Sigma \setminus
\text{Int}(N(\alpha))$. Thus $H_2(\Sigma, \alpha) = 0$ if and only
if for every such component $C$ the group $H_2(C,\partial
N(\alpha)\cap C) = 0$, ie, when $C \cap \partial \Sigma \neq
\emptyset$. Thus $H_2(\Sigma, \alpha) =0$ exactly when the map
$\pi_0(\partial \Sigma) \to \pi_0(\Sigma \setminus \alpha)$ is
injective.
\end{proof}

\fullref{prop:1} justifies the following definition.

\begin{defn} \label{defn:5}
A sutured Heegaard diagram $( \Sigma, \boldsymbol{\alpha},
\boldsymbol{\beta})$ is called \emph{balanced\/} if
$|\boldsymbol{\alpha}| = |\boldsymbol{\beta}|$ and the maps
$\pi_0(\partial \Sigma) \to \pi_0(\Sigma \setminus \bigcup
\boldsymbol{\alpha})$ and $\pi_0(\partial \Sigma) \to \pi_0(\Sigma
\setminus \bigcup \boldsymbol{\beta})$ are surjective.
\end{defn}

\begin{rem}
We will use the abbreviation ``\emph{balanced diagram\/}'' for ``balanced
sutured Heegaard diagram''.
\end{rem}

\begin{prop} \label{prop:2}
Let $(M,\gamma)$ be a sutured manifold for which the maps
$$\pi_0(R_+(\gamma)) \to \pi_0(M)\,\,\, \text{and}\,\,\, \pi_0(R_-(\gamma)) \to
\pi_0(M)$$ are surjective. Then there exists a sutured Heegaard
diagram $(\Sigma, \boldsymbol{\alpha}, \boldsymbol{\beta})$ defining
it.
\end{prop}

\begin{proof}
Fix a Riemannian metric on $M$. First we construct a special Morse
function $f$ on $M$ the following way. Choose a diffeomorphism
$\varphi \co \gamma \to s(\gamma) \times [-1,4]$ so that
$\varphi(s(\gamma)) = s(\gamma) \times \{3/2\}$ and let $p_2 \co
s(\gamma) \times [-1,4] \to [-1,4]$ be the projection onto the
second factor. Then we define $f|\gamma$ to be $p_2 \circ \varphi$. 
Furthermore, let $f|R_-(\gamma) \equiv -1$ and $f|R_+(\gamma) \equiv
4$. Now take a generic extension of $f|\partial M$ to $M$. Then $f
\co M \to \mathbb{R}$ is a Morse function.

Using \cite[Theorem 4.8]{Milnor} we can assume that $f$ is
self-indexing. Applying the idea of \cite[Theorem 8.1]{Milnor} as
follows we can assume that $f$ has no index $0$ and $3$ critical
points. Since the map $H_0(R_-(\gamma)) \to H_0(M)$ is surjective
$H_0(M, R_-(\gamma)) =0$. Thus, using CW homology, we see that for
every index $0$ critical point of $f$ we can find an index $1$
critical point so that there is exactly one gradient flow line
connecting them, and they can be canceled. Indeed, since $H_0(M,
R_-(\gamma)) =0$, for every index zero critical point $p$ there is
an index one critical point $q$ such that $p$ and $q$ are connected
by an odd number of gradient flow lines. But there are only two flow
lines coming out of $q$, so there is exactly one trajectory
connecting $p$ and $q$. During this process we do not have to change
$f|\partial M$. Similarly, we can cancel every index $3$ critical
point of $f$. 

Finally, let $\Sigma = f^{-1}(3/2)$ and let $\boldsymbol{\alpha}$
and $\boldsymbol{\beta}$ be the intersections of $\Sigma$ with the
ascending and descending manifolds of the index one and two critical
points of $f$ respectively. Then $(\Sigma, \boldsymbol{\alpha},
\boldsymbol{\beta})$ defines $(M,\gamma)$. 
\end{proof}

\begin{prop} \label{prop:3}
For every balanced sutured manifold $(M,\gamma)$ there exists a
balanced diagram defining it.
\end{prop}

\begin{proof}
This is a corollary of \fullref{prop:1} and \fullref{prop:2}.
\end{proof}

Next we will state and prove a generalization of \cite[Proposition
2.2]{OSz}.

\begin{prop} \label{prop:4}
If the balanced diagrams $(\Sigma_0, \boldsymbol{\alpha_0},
\boldsymbol{\beta_0})$ and $(\Sigma_1, \boldsymbol{\alpha_1},
\boldsymbol{\beta_1})$ define the same balanced sutured manifold
$(M, \gamma)$ then they are diffeomorphic after a finite sequence of
Heegaard moves.
\end{prop}

\begin{proof}
Suppose that $\alpha \subset R_-(\gamma)$ is a simple closed
curve such that the $1$--handle attached to $M$ along $\alpha$ can be
canceled by a $0$--handle $B^3$. Then the curve $\alpha$ bounds the
$2$--disc $\partial B^3 \cap R_-(\gamma)$. 

Using the above observation we get that adding a canceling pair of
index 0 and 1 critical points corresponds to adding a curve $\alpha$
to $\boldsymbol{\alpha}$ such that after performing surgery on
$\Sigma$ along $\boldsymbol{\alpha}$ (so that we obtain
$R_-(\gamma)$) the image of $\alpha$ bounds a disc.

\begin{note}
If $\boldsymbol{\gamma}$ is a set of pairwise disjoint simple closed
curves in the interior of a surface $\Sigma$ then
$\Sigma[\boldsymbol{\gamma}]$ denotes the surface obtained by
surgery on $\Sigma$ along $\boldsymbol{\gamma}$. 
\end{note}

\begin{lem} \label{lem:2}
Let $\alpha_1, \dots, \alpha_d, \gamma$ and $\delta$ be pairwise
disjoint simple closed curves in a compact oriented surface $\Sigma$
such that the image of both $\gamma$ and $\delta$ bound a disc in
$\Sigma[\alpha_1, \dots, \alpha_d]$. Suppose that $\gamma$ is not
null-homologous.

Then  there is an $i \in \{\,1, \dots, d\,\}$ such that $\gamma$ is isotopic to
a curve obtained by handlesliding $\alpha_i$ across some collection
of the $\alpha_j$ for $j \neq i$. Moreover, the curves $\alpha_i$
and $\delta$ both bound discs in $\widetilde{\Sigma} =
\Sigma[\alpha_1, \dots, \alpha_{i-1}, \alpha_{i+1}, \dots, \alpha_d,
\gamma]$. 
\end{lem}

\begin{proof}
Let $D_{\gamma}$ and $D_{\delta}$ be discs bound by $\gamma$ and
$\delta$ in $\Sigma' = \Sigma[\alpha_1, \dots, \alpha_d]$
respectively. For $i=1, \dots, d$ let $p_i, q_i \in \Sigma'$ be the
points corresponding to the zero-sphere which replaced the circle
$\alpha_i$. Since $\gamma$ is not null-homologous, there is an $i
\in \{\,1, \dots, d\,\}$ such that $D_{\gamma}$ separates $p_i$ and
$q_i$. We can suppose without loss of generality that $i = 1$ and
$p_1 \in D_{\gamma}$ while $q_1 \not\in D_{\gamma}$. An isotopy in
$D_{\gamma}$ of a small circle around $p_1$ to $\gamma$ corresponds
to handlesliding $\alpha_1$ across some collection of the $\alpha_j$
for $j \neq 1$ so that we obtain $\gamma$. 

Observe that $\widetilde{\Sigma}$ is obtained from $\Sigma'$ by
adding a tube $T$ to $\Sigma' \setminus \{\, p_1, q_1\,\}$ and
performing surgery along $\gamma$. We take $\Sigma' \setminus
\gamma$ and pinch the boundary component corresponding to $\partial
D_{\gamma}$ to $p_0$ and $\partial (\Sigma' \setminus D_{\gamma})$
to $q_0$. Then $\alpha_1$ is the boundary of the disc $(D_{\gamma}
\setminus \{p_1\}) \cup \{p_0\} \subset \widetilde{\Sigma}$. 

We are now going to prove that $\delta$ bounds a disc in
$\widetilde{\Sigma}$. If $p_1 \not\in D_{\delta}$ and $q_1 \not\in
D_{\delta}$ then since $\delta \cap \gamma = \emptyset$ the disc
$D_{\delta}$ ``survives'' in $\widetilde{\Sigma}$. If $p_1 \in
D_{\delta}$ and $q_1 \not \in D_{\delta}$ then in
$\widetilde{\Sigma}$ the curve $\delta$ bounds $(D_{\gamma}
\setminus D_{\delta}) \cup \{p_0\}$ if $D_{\gamma} \supset
D_{\delta}$ and $(D_{\delta} \setminus D_{\gamma}) \cup \{q_0\}$
otherwise. If $p_1 \not\in D_{\delta}$ and $q_1 \in D_{\delta}$ then
$\delta$ bounds in $\widetilde{\Sigma}$ the disc $(D_{\delta}
\setminus \{q_1\}) \cup T \cup (D_{\gamma} \setminus \{p_1\}) \cup
\{p_0\}$. Finally, if $p_1, q_1 \in D_{\delta}$ then of course
$D_{\gamma} \subset D_{\delta}$, and in $\widetilde{\Sigma}$ the
curve $\delta$ bounds $(D_{\delta} \setminus (D_{\gamma} \cup
\{q_1\})) \cup T \cup (D_\gamma \setminus \{p_1\})\cup \{p_0,
q_0\}$. 

(In fact, $\Sigma' = \Sigma[\alpha_2, \dots, \alpha_d][\alpha_1]$
and $\widetilde{\Sigma} = \Sigma[\alpha_2, \dots,
\alpha_d][\gamma]$, and furthermore the curves $\alpha_1$ and $\gamma$
are isotopic in $\Sigma[\alpha_2, \dots, \alpha_d]$.)
\end{proof}

\begin{lem} \label{lem:3}
Let $\boldsymbol{\delta}$ be a set of pairwise disjoint simple
closed curves in $\Sigma$, and suppose that we are given two subsets
of curves $\boldsymbol{\alpha}, \boldsymbol{\gamma} \subset
\boldsymbol{\delta}$ that are linearly independent in $H_1(\Sigma;
\mathbb{Q})$. Suppose furthermore that the image of every $\delta
\in \boldsymbol{\delta} \setminus \boldsymbol{\alpha}$ bounds a disc
in $\Sigma[\boldsymbol{\alpha}]$. Then $\boldsymbol{\gamma}$ can be
obtained from $\boldsymbol{\alpha}$ by a series of isotopies and
handleslides. Moreover, the image of every $\delta \in
\boldsymbol{\delta} \setminus \boldsymbol{\gamma}$ bounds a disc in
$\Sigma[\boldsymbol{\gamma}]$. 
\end{lem}

\begin{proof}
Let $d = |\boldsymbol{\alpha}| = |\boldsymbol{\gamma}|$. We prove
the claim using induction on $d$. The case $d=0$ is trivial. Note
that it follows from the hypothesis that $\boldsymbol{\alpha}$ and
$\boldsymbol{\gamma}$ span the same subspace in $H_1(\Sigma;
\mathbb{Q})$. 

If $\boldsymbol{\alpha} \cap \boldsymbol{\gamma} \neq \emptyset$, 
say the curve $\alpha$ lies in the intersection, then perform
surgery on $\Sigma$ along $\alpha$ to obtain a new surface $\Sigma'$
with two marked points $p, q$ and two $(d-1)$--tuples of curves
$\boldsymbol{\alpha}'$ and $\boldsymbol{\gamma}'$. Let
$\boldsymbol{\delta}' = \boldsymbol{\delta} \setminus \{\alpha\}$. 
Note that every $\delta \in \boldsymbol{\delta'} \setminus
\boldsymbol{\alpha'} = \boldsymbol{\delta} \setminus
\boldsymbol{\alpha}$ bounds a disc in $\Sigma'[\boldsymbol{\alpha'}]
= \Sigma[\boldsymbol{\alpha}]$. Using the induction hypothesis
$\boldsymbol{\alpha}'$ and $\boldsymbol{\gamma}'$ are related by
isotopies and handleslides. We can arrange (using isotopies) that
each handleslide is disjoint from $p$ and $q$. Each isotopy of a
curve in $\Sigma'$ that crosses $p$ or $q$ corresponds to a
handleslide in $\Sigma$ across $\alpha$. Thus $\boldsymbol{\alpha}$
and $\boldsymbol{\gamma}$ are also related by isotopies and
handleslides. Also from the induction hypothesis we get that every
$\delta \in \boldsymbol{\delta'} \setminus \boldsymbol{\gamma'}$
bounds a disc in $\Sigma'[\boldsymbol{\gamma'}] =
\Sigma[\boldsymbol{\gamma}]$. This implies that the image of every
$\delta \in \boldsymbol{\delta} \setminus \boldsymbol{\gamma}$
bounds a disc in $\Sigma[\boldsymbol{\gamma}]$. 

If $\boldsymbol{\alpha} \cap \boldsymbol{\gamma} = \emptyset$ then
take any $\gamma \in \boldsymbol{\gamma}$. Since elements of
$\boldsymbol{\gamma}$ are linearly independent $\gamma$ is not
null-homologous. Thus, using \fullref{lem:2}, $\gamma$ can be
obtained by handlesliding some $\alpha_i$ across a collection of the
$\alpha_j$ for $j \neq i$. So we have reduced to the case where the
two subsets are not disjoint.
\end{proof}

For $i \in \{\,0,1\,\}$ choose a Morse function $f_i$ inducing
$(\Sigma_i, \boldsymbol{\alpha_i}, \boldsymbol{\beta_i})$ as in the
proof of \fullref{prop:2} and let $\{\,f_t \,\colon\, 0 \le t
\le 1\,\}$ be a generic one-parameter family of functions connecting
them. We can suppose that $f_t$ is fixed in a neighborhood of
$\partial M$. Also equip $Y$ with a generic Riemannian metric. Then
there is a finite subset $E \subset I$ such that for $t \in I
\setminus E$ the function $f_t$ is Morse with gradient flow lines
flowing only from larger to strictly smaller index critical points,
and thus induces a diagram $(\Sigma_t, \boldsymbol{\alpha_t},
\boldsymbol{\beta_t})$. Here $\boldsymbol{\alpha_t}$ and
$\boldsymbol{\beta_t}$ are the intersections of $\Sigma_t$ with the
ascending and descending manifolds of the index one and two critical
points of $f_t$ respectively. As $t$ passes through an element $e
\in E$ the diagram corresponding to $f_t$ experiences one of the
following changes. There is either a handleslide among the $\alpha$
curves or the $\beta$ curves (corresponding to a gradient flow line
connecting two index one or two index two critical points of $f_e$),
or a stabilization/destabilization (corresponding to
creation/cancellation of index 1 and 2 critical points), or a new
$\alpha$ or $\beta$ curve appears/disappears (corresponding to
canceling index 0 and 1, or index 2 and 3 critical points). The last
case is called a pair creation/cancellation.

For each $t \in I \setminus E$ choose two maximal homologically
linearly independent subsets $\boldsymbol{\alpha_t'} \subset
\boldsymbol{\alpha_t}$ and $\boldsymbol{\beta_t'} \subset
\boldsymbol{\beta_t}$ that change continuously in $t$. Then of
course $\boldsymbol{\alpha_i'} = \boldsymbol{\alpha_i}$ and
$\boldsymbol{\beta_i'} = \boldsymbol{\beta_i}$ for $i=0,1$. So it is
enough to show that for every $e \in E$ and sufficiently small
$\varepsilon$ the diagrams $(\Sigma_{e-\varepsilon},
\boldsymbol{\alpha_{e-\varepsilon}'},
\boldsymbol{\beta_{e-\varepsilon}'})$ and $(\Sigma_{e+\varepsilon},
\boldsymbol{\alpha_{e+\varepsilon}'},
\boldsymbol{\beta_{e+\varepsilon}'})$ are related by isotopies,
handleslides, stabilization and destabilization.

In order to do this we also need to prove the fact that for every $t
\in I \setminus E$ and every curve $\alpha \in \boldsymbol{\alpha_t}
\setminus \boldsymbol{\alpha_t'}$ the image of $\alpha$ bounds a
disc in the surface $\Sigma_t[\boldsymbol{\alpha_t'}]$. We prove
this by induction on the component of $I \setminus E$ containing
$t$. It is obviously true for $t=0$. 

First consider the case when $e$ does not correspond to stabilization
or destabilization. Let $\boldsymbol{\delta} =
\boldsymbol{\alpha_{e-\varepsilon}} \cup
\boldsymbol{\alpha_{e+\varepsilon}}$, this is a set of pairwise
disjoint curves. Furthermore, let $\boldsymbol{\alpha} =
\boldsymbol{\alpha_{e-\varepsilon}'}$ and $\boldsymbol{\gamma} =
\boldsymbol{\alpha_{e+\varepsilon}'}$. Observe that
$\boldsymbol{\delta} \setminus \boldsymbol{\alpha_{e-\varepsilon}}$
consists of at most one curve $\delta$ obtained from either a
handleslide within $\boldsymbol{\alpha_{e-\varepsilon}}$ or a pair
creation. Using the induction hypothesis for $t = e-\varepsilon$ we
see that $\delta$ also bounds a disc in
$\Sigma_{e-\varepsilon}[\boldsymbol{\alpha}]$. Thus we can apply
\fullref{lem:3} to $\boldsymbol{\alpha}, \boldsymbol{\gamma}
\subset \boldsymbol{\delta}$ showing that $\boldsymbol{\alpha}$ and
$\boldsymbol{\gamma}$ are related by isotopies and handleslides and
that the induction hypothesis also holds for $t = e + \varepsilon$. 
A similar argument applies to the $\beta$ curves.

Now suppose that $e$ corresponds to a stabilization; the new
curves appearing are $\alpha$ and $\beta$. Define
$\boldsymbol{\alpha} = \boldsymbol{\alpha_{e- \varepsilon}'} \cup
\{\alpha\}$ and $\boldsymbol{\beta} = \boldsymbol{\beta_{e-
\varepsilon}'} \cup \{\beta\}$, considered as sets of curves in
$\Sigma_{e+\varepsilon}$. Then we can apply \fullref{lem:3} to
$\boldsymbol{\alpha}, \boldsymbol{\alpha'_{e+\varepsilon}} \subset
\boldsymbol{\alpha_{e+\varepsilon}}$ and $\boldsymbol{\beta},
\boldsymbol{\beta'_{e+\varepsilon}} \subset
\boldsymbol{\beta_{e+\varepsilon}}$. 

The case of a destabilization is proved in a similar way, by taking
$\boldsymbol{\alpha} = \boldsymbol{\alpha_{e+ \varepsilon}'} \cup
\{\alpha\}$ and $\boldsymbol{\beta} = \boldsymbol{\beta_{e+
\varepsilon}'} \cup \{\beta\}$, where $\alpha$ and $\beta$ are the
curves that vanish.
\end{proof}

\begin{rem} \label{rem:1}
From \fullref{prop:4} we see that if we associate to every
balanced diagram a quantity that is unchanged by isotopies,
handleslides and stabilization we get a topological invariant of
sutured $3$--manifolds.
\end{rem}

\section{Whitney discs and their domains}

For a surface $\Sigma$ let $\text{Sym}^d(\Sigma)$ denote the d-fold
symmetric product $\Sigma^{\times d} / S_d$. This is always a smooth
$2d$--manifold. A complex structure $\mathfrak{j}$ on $\Sigma$
naturally endows $\text{Sym}^d(\Sigma)$ with a complex structure,
denoted $\text{Sym}^d(\mathfrak{j})$. This structure
$\text{Sym}^d(\mathfrak{j})$ is specified by the property that the
quotient map $\Sigma^d \to \text{Sym}^d(\Sigma)$ is holomorphic.

\begin{defn} \label{defn:6}
Let $(\Sigma, \boldsymbol{\alpha}, \boldsymbol{\beta})$ be a
balanced diagram, where $\boldsymbol{\alpha} = \{\, \alpha_1, \dots,
\alpha_d \,\}$ and $\boldsymbol{\beta} = \{\, \beta_1, \dots,
\beta_d \,\}$. Then let $\mathbb{T}_{\alpha} = (\alpha_1 \times
\dots \times \alpha_d) / S_d$ and $\mathbb{T}_{\beta} = (\beta_1
\times \dots \times \beta_d) / S_d$. These are $d$--dimensional tori
in $\text{Sym}^{\smash{d}}(\Sigma)$. 
\end{defn}

\begin{lem} \label{lem:4}
For a balanced diagram $(\Sigma, \boldsymbol{\alpha},
\boldsymbol{\beta})$ and an arbitrary complex structure
$\mathfrak{j}$ on $\Sigma$, the submanifolds $\mathbb{T}_{\alpha},
\mathbb{T}_{\beta} \subset \text{Sym}^d(\Sigma)$ are totally real
with respect to $\text{Sym}^d(\mathfrak{j})$. 
\end{lem}

\begin{proof}
The submanifold $\alpha_1 \times \dots \times \alpha_d \subset
\Sigma^{\times d}$ is totally real with respect to
$\mathfrak{j}^{\times d}$ and misses the diagonal (consisting of
those $d$--tuples for which at least two coordinates coincide). The
claim thus follows since the projection map $\pi \co
\Sigma^{\times d} \to \text{Sym}^d(\Sigma)$ is a holomorphic local
diffeomorphism away from the diagonal.
\end{proof}

Note that if every $\alpha \in \boldsymbol{\alpha}$ and $\beta \in
\boldsymbol{\beta}$ are transversal then the tori
$\mathbb{T}_{\alpha}$ and $\mathbb{T}_{\beta}$ intersect
transversally.

\begin{note} \label{note:1}
Let $\mathbb{D}$ denote the unit disc in $\mathbb{C}$, and let
$e_1=\{\, z \in \partial \mathbb{D} \,\colon\, \text{Re}(z) \ge 0\,\}$
and $e_2=\{\, z \in \partial \mathbb{D} \,\colon\, \text{Re}(z) \le
0\,\}$. 
\end{note}

\begin{defn} \label{defn:7}
Let $\mathbf{x}, \mathbf{y} \in \mathbb{T}_{\alpha} \cap
\mathbb{T}_{\beta}$ be intersection points. A \emph{Whitney disc
connecting $\mathbf{x}$ to $\mathbf{y}$} is a continuous map $u
\co \mathbb{D} \to \text{Sym}^{\smash{d}}(\Sigma)$ such that $u(-i) =
\mathbf{x}$, $u(i) = \mathbf{y}$ and $u(e_1) \subset
\mathbb{T}_{\alpha}$, $u(e_2) \subset \mathbb{T}_{\beta}$. Let
$\pi_2(\mathbf{x}, \mathbf{y})$ denote the set of homotopy classes
of Whitney discs connecting $\mathbf{x}$ to $\mathbf{y}$. 
\end{defn}

\begin{defn} \label{defn:8}
For $z \in \Sigma \setminus (\bigcup \boldsymbol{\alpha} \cup \bigcup
\boldsymbol{\beta})$ and $u$ a Whitney disc, choose a
Whitney disc $u'$ homotopic to $u$ such that $u'$ intersects the
hypersurface $ \{z\} \times \text{Sym}^{\smash{d-1}}(\Sigma)$ transversally.
Define $n_z(u)$ to be the algebraic intersection number $u' \cap
(\{z\} \times \text{Sym}^{d-1}(\Sigma))$. 
\end{defn}

Note that $n_z(u)$ only depends on the component of $\Sigma
\setminus (\bigcup \boldsymbol{\alpha} \cup \bigcup
\boldsymbol{\beta})$ in which $z$ lies and on the homotopy class of
the Whitney disc $u$. Moreover, if the component of $z$ contains a
boundary component of $\Sigma$ then $n_z(u) = 0$. Indeed, we can
choose $z$ on $\partial \Sigma$ and we can homotope $u$ to be
disjoint from $\partial \text{Sym}^d(\Sigma) \supset \{z\} \times
\text{Sym}^{d-1}(\Sigma)$, showing that $n_z(u) = 0$. This last
remark implies that we can run the Floer homology machinery without
worrying about being in a manifold with boundary.

\begin{defn} \label{defn:9}
For a balanced diagram let $\mathcal{D}_1, \dots, \mathcal{D}_m$
denote the closures of the components of $\Sigma \setminus (\bigcup
\boldsymbol{\alpha} \cup \bigcup \boldsymbol{\beta})$ disjoint from
$\partial \Sigma$. Then let $D(\Sigma, \boldsymbol{\alpha},
\boldsymbol{\beta})$ be the free abelian group generated by
$\{\,\mathcal{D}_1, \dots, \mathcal{D}_m\,\}$. This is of course
isomorphic to $\mathbb{Z}^m$. We call an element of $D(\Sigma,
\boldsymbol{\alpha}, \boldsymbol{\beta})$ a \emph{domain\/}. An
element $\mathcal{D}$ of $\mathbb{Z}_{\ge 0}^m$ is called a
\emph{positive\/} domain, we write $\mathcal{D} \ge 0$. A domain
$\mathcal{P} \in D(\Sigma, \boldsymbol{\alpha}, \boldsymbol{\beta})$
is called a \emph{periodic domain\/} if the boundary of the $2$--chain
$\mathcal{P}$ is a sum of $\alpha$- and $\beta$--curves.
\end{defn}

\begin{defn} \label{defn:10}
For every $1 \le i \le m$ choose a point $z_i \in \mathcal{D}_i$. 
Then the \emph{domain of a Whitney disc $u$\/} is defined as
$$\mathcal{D}(u) = \sum_{i=1}^m n_{z_i}(u)\mathcal{D}_i \in
D(\Sigma, \boldsymbol{\alpha}, \boldsymbol{\beta}).$$ For $\phi \in
\pi_2(\mathbf{x}, \mathbf{y})$ and $u$ a representative of the
homotopy class $\phi$, let $\mathcal{D}(\phi) = \mathcal{D}(u)$. 
\end{defn}

\begin{rem} \label{rem:2}
If a Whitney disc $u$ is holomorphic then $\mathcal{D}(u) \ge 0$. 

If $\mathbf{x} \in \mathbb{T}_{\alpha} \cap \mathbb{T}_{\beta}$ and
if $u$ is a Whitney disc connecting $\mathbf{x}$ to $\mathbf{x}$
then $\mathcal{D}(u)$ is a periodic domain.
\end{rem}

\begin{defn} \label{defn:15}
If $(\Sigma, \boldsymbol{\alpha}, \boldsymbol{\beta})$ is a balanced
diagram defining the balanced sutured manifold $(M, \gamma)$ and if
$\mathcal{P} \in D(\Sigma, \boldsymbol{\alpha}, \boldsymbol{\beta})$
is a periodic domain then we can naturally associate to
$\mathcal{P}$ a homology class $H(\mathcal{P}) \in H_2(M;
\mathbb{Z})$ as follows. The boundary of the two-chain $\mathcal{P}$
is a sum $\sum_{i=1}^d a_i\alpha_i + \sum_{i=1}^d b_i \beta_i$. Let
$A_i$ denote the core of the two-handle attached to $\alpha_i$ and
$B_i$ the core of the two-handle attached to $\beta_i$. Then let
$$H(\mathcal{P}) = \left[\mathcal{P} + \sum_{i=1}^d a_iA_i +
\sum_{i=1}^d b_iB_i\right] \in H_2(M; \mathbb{Z}).$$
\end{defn}

\begin{lem} \label{lem:6}
If $H(\mathcal{P}) =0$ then $\mathcal{P} = 0$. 
\end{lem}

\begin{proof}
Since $\Sigma$ has no closed components we have that $H_2(\Sigma;
\mathbb{Z}) = 0$. Thus, if $\mathcal{P} \neq 0$ then $\partial
\mathcal{P} = \sum_{i=1}^d a_i\alpha_i + \sum_{i=1}^d b_i \beta_i
\neq 0$. Suppose for example that $a_1 \neq 0$. This implies
that $H(\mathcal{P})$ has nonzero algebraic intersection with the
co-core $A_1'$ of the two-handle attached to $\alpha_1$ (whose core
is $A_1$). Since $[A_1'] \neq 0$ in $H_1(M, \partial M; \mathbb{Z})$
we get that $H(\mathcal{P}) \neq 0$. 
\end{proof}

\begin{defn} \label{defn:17}
A balanced diagram $(\Sigma, \boldsymbol{\alpha},
\boldsymbol{\beta})$ is called \emph{admissible\/} if every periodic
domain $\mathcal{P} \neq 0$ has both positive and negative
coefficients.
\end{defn}

\begin{cor} \label{cor:1}
If $(M, \gamma)$ is a balanced sutured manifold such that $H_2(M;
\mathbb{Z})$ is $0$ and if $(\Sigma, \boldsymbol{\alpha},
\boldsymbol{\beta})$ is an arbitrary balanced diagram defining
$(M,\gamma)$ then there are no nonzero periodic domains in
$D(\Sigma, \boldsymbol{\alpha, \boldsymbol{\beta}})$. Thus any
balanced diagram defining $(M, \gamma)$ is automatically admissible.
\end{cor}

\begin{defn} \label{defn:16}
Let $\mathbf{x}, \mathbf{y} \in \mathbb{T}_{\alpha} \cap
\mathbb{T}_{\beta}$. A domain $\mathcal{D} \in D(\Sigma,
\boldsymbol{\alpha}, \boldsymbol{\beta})$ is said to \emph{connect
$\mathbf{x}$ to $\mathbf{y}$} if for every $1 \le i \le d$ the
equalities $\partial(\alpha_i \cap \partial \mathcal{D}) =
(\mathbf{x} \cap \alpha_i) - (\mathbf{y} \cap \alpha_i)$ and
$\partial(\beta_i \cap \partial \mathcal{D}) = (\mathbf{x} \cap
\beta_i) - (\mathbf{y} \cap \beta_i)$ hold. We will denote
by $D(\mathbf{x}, \mathbf{y})$ the set of domains connecting
$\mathbf{x}$ to $\mathbf{y}$. 
\end{defn}

Note that if $\phi \in \pi_2(\mathbf{x}, \mathbf{y})$ then
$\mathcal{D}(\phi) \in D(\mathbf{x}, \mathbf{y})$. 

\begin{lem} \label{lem:7}
If the balanced diagram $(\Sigma, \boldsymbol{\alpha},
\boldsymbol{\beta})$ is admissible then for every pair $\mathbf{x},
\mathbf{y} \in \mathbb{T}_{\alpha} \cap \mathbb{T}_{\beta}$ the set
$\{\, \mathcal{D} \in D(\mathbf{x}, \mathbf{y}) \,\colon\,
\mathcal{D} \ge 0 \,\}$ is finite.
\end{lem}

\begin{proof}
The argument that follows can be found in the proof of \cite[Lemma
4.13]{OSz}. If $D(\mathbf{x}, \mathbf{y}) \neq \emptyset$ then fix
an element $\mathcal{D}_0 \in D(\mathbf{x}, \mathbf{y})$. Then every
element $\mathcal{D} \in D(\mathbf{x}, \mathbf{y})$ can be written
as $\mathcal{D} = \mathcal{D}_0 + \mathcal{P}$, where $\mathcal{P}
\in D(\mathbf{x}, \mathbf{x})$ is a periodic domain. Hence if
$\mathcal{D} \ge 0$ then $\mathcal{P} \ge -\mathcal{D}_0$. 

So the lemma follows if we show that the set $Q = \{\,\mathcal{P}
\in D(\mathbf{x}, \mathbf{x}) \,:\, \mathcal{P} \ge
-\mathcal{D}_0 \,\}$ is finite. We can think of $Q$ as a subset of
the lattice $\mathbb{Z}^m \subset \mathbb{R}^m$. If $Q$ had
infinitely many elements, then we could find a sequence
$(p_j)_{j=1}^{\infty}$ in $Q$ with $\|p_j\| \to \infty$. Taking a
subsequence we can suppose that $(p_j/\|p_j\|)$ converges to a unit
vector $p$ in the vector space of periodic domains with real
coefficients. Since the coefficients of $p_j$ are bounded below and
$\|p_j\| \to \infty$ we get that $p \ge 0$. Thus the polytope
consisting of vectors corresponding to real periodic domains with
$\ge 0$ multiplicities also has a nonzero rational vector. After
clearing denominators we obtain a nonzero integer periodic domain
with nonnegative multiplicities. This contradicts the hypothesis of
admissibility.
\end{proof}

\begin{prop} \label{prop:5}
Every balanced diagram $(\Sigma, \boldsymbol{\alpha},
\boldsymbol{\beta})$ is isotopic to an admissible one.
\end{prop}

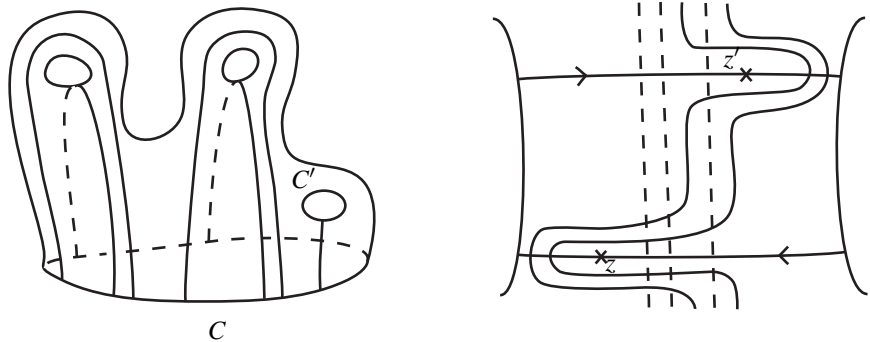
\begin{figure}[ht!]
\labellist
\small\hair 2pt
\pinlabel $z'$ [br] at 420 519
\pinlabel $z$ [tl] at 366 450
\pinlabel $C$ [t] at 221 427
\pinlabel $C'$ [br] at 260 474
\endlabellist
\centering
\includegraphics{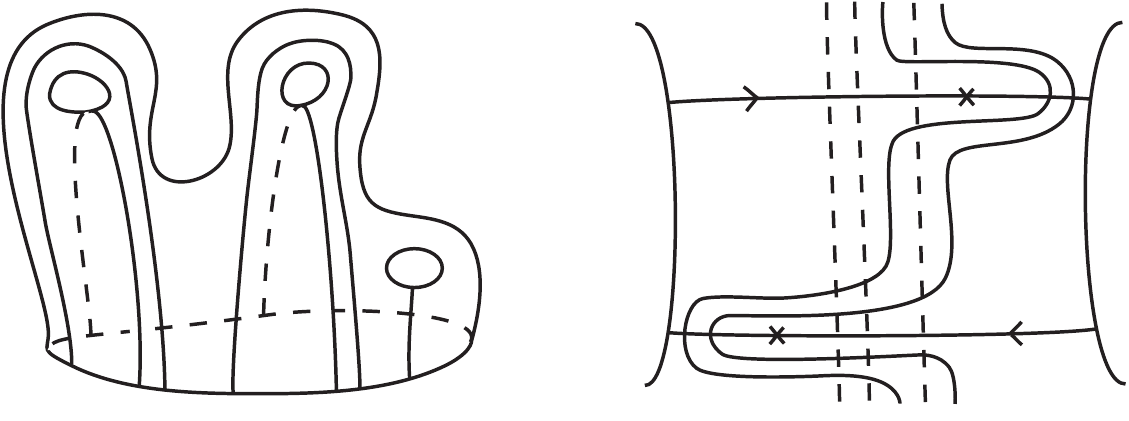}
\caption{The picture on the left shows curves that represent a basis
of $H_1(\Sigma,\partial\Sigma; \mathbb{Z})$. On the right we can see
the procedure to achieve\break admissibility. The dotted lines represent
the $\alpha$ curves, and the solid lines represent the $\beta$
curves.} \label{fig:1}
\end{figure}

\begin{proof}
Fix a boundary component $C \subset \partial \Sigma$. We can choose
a set of pairwise disjoint, oriented and properly embedded arcs
$\gamma_1, \dots, \gamma_l \subset \Sigma$ such that for every $1
\le i \le l$ the endpoints $\partial \gamma_i$ lie in
$\partial\Sigma$; furthermore these arcs generate the relative
homology group $H_1(\Sigma, \partial\Sigma; \mathbb{Z})$. This can
be done as follows (see the left hand side of \fullref{fig:1}).
Let $\Sigma'$ denote the surface obtained from $\Sigma$ by gluing a
disc to every component of $\partial \Sigma \setminus C$. Let $g$
denote the genus of $\Sigma'$. Then we can choose a set of $2g$
curves in $\Sigma'$ as above, that are also disjoint from
$\partial\Sigma \setminus C$. Finally, for each component $C'$ of
$\partial \Sigma \setminus C$ connect $C$ and $C'$ with a $\gamma$
curve. Note that $\Sigma \setminus (\partial \Sigma \cup
\bigcup_{i=1}^l \gamma_i)$ is homeomorphic to an open disc.

We perform an isotopy of the $\beta$ curves in a regular
neighborhood of $\gamma_1 \cup \dots \cup \gamma_l$ as described in
the proof of \cite[Proposition 3.6]{OSz2}. Specifically, for every
$1 \le i \le l$ choose an oppositely oriented parallel copy
$\gamma_i'$ of $\gamma_i$. Using a finger move isotope the $\beta$
curves intersecting $\gamma_i$ towards the endpoint of $\gamma_i$ so
that there is a point $z_i \in \gamma_i$ separating
$\boldsymbol{\alpha} \cap \gamma_i$ from $\boldsymbol{\beta} \cap
\gamma_i$. Perform a similar isotopy of the $\beta$ curves in a
neighborhood of each $\gamma_i'$. The point $z_i' \in \gamma_i'$
separates $\boldsymbol{\alpha} \cap \gamma_i'$ from
$\boldsymbol{\beta} \cap \gamma_i'$. See the right-hand side of
\fullref{fig:1}.

We claim the diagram obtained this way is admissible. Let
$\mathcal{P}$ be a periodic domain. Then $$\partial \mathcal{P} =
\sum_{i=1}^d a_i \cdot \alpha_i + \sum_{i=1}^d b_i \cdot \beta_i = A
+ B.$$ First suppose that there is an $1 \le i \le l$ such that the
algebraic intersection $A \cap \gamma_i \neq 0$. Since the
multiplicity of $\mathcal{P}$ at the points of $\partial \Sigma$ is
$0$ we get that the multiplicity of $\mathcal{P}$ at $z_i$ is $A
\cap \gamma_i$ and at $z_i'$ it is $A \cap \gamma_i' = -A \cap
\gamma_i$ (see the right-hand side of \fullref{fig:1}). Indeed,
$z_i$ separates $A \cap \gamma_i$ from $B \cap \gamma_i$ on
$\gamma_i$ and $z_i'$ separates $A \cap \gamma_i'$ from $B \cap
\gamma_i'$ on $\gamma_i'$. Thus $\mathcal{P}$ has both positive and
negative multiplicities.

On the other hand, if for every $1 \le i \le l$ the intersection
number $A \cap \gamma_i =0$, then since $\gamma_1, \dots, \gamma_l$
span $H_1(\Sigma, \partial \Sigma; \mathbb{Z})$, we get that $A$ is
null-homologous in $\Sigma$. Indeed, in this case $A$ is homologous
to a curve lying in $\Sigma \setminus (\gamma_1 \cup \dots \cup
\gamma_l) \approx D^2$. Since the elements of $\boldsymbol{\alpha}$
are linearly independent is $H_1(\Sigma; \mathbb{Z})$ we get that
for every $1 \le j \le d$ the coefficient $a_j = 0$. But $\partial
\mathcal{P} = A+B$ implies that $B$ is homologous with $-A$ in
$\Sigma$, thus $B \sim 0$. So we get that $b_j = 0$ for every $1 \le
j \le d$. Thus in this case $\mathcal{P} = 0$. 
\end{proof}

\section[Spin-c structures]{$\text{Spin}^c$ structures}

In this section $(M,\gamma)$ denotes a connected balanced sutured
manifold.

\begin{note} \label{note:2}
Let $v_0$ be the nonzero vector field along $\partial M$ that
points into $M$ along $R_-(\gamma)$, points out of $M$ along
$R_+(\gamma)$, and on $\gamma$ it is the gradient of the height
function $s(\gamma) \times I \to I$. The space of such vector fields
is contractible.

The field $v_0^{\perp}$ is an oriented two-plane field along $\partial M$. We
will use the notation $$\delta = c_1(v_0^{\perp}) = e(v_0^{\perp})
\in H^2(\partial M; \mathbb{Z}).$$
\end{note}

\begin{defn} \label{defn:11}
Let $v$ and $w$ be vector fields on $M$ that agree with $v_0$ on
$\partial M$. We say that $v$ and $w$ are \emph{homologous\/} if there
is an open ball $B \subset \text{Int}(M)$ such that $v|(M \setminus
B)$ is homotopic to $w|(M \setminus B)$ rel $\partial M$. We define
$\text{Spin}^c(M, \gamma)$ to be the set of homology classes of
nonzero vector fields $v$ on $M$ such that $v|\partial M = v_0$. 
\end{defn}

\begin{rem} \label{rem:3}
Let $f$ be a Morse function as in \fullref{prop:2}. Then the
vector field $\text{grad}(f)|\partial M = v_0$, the number $d$ of
index 1 and 2 critical points of $f$ agree, and $f$ has no index 0
or 3 critical points. Choose $d$ pairwise disjoint balls in $M$, 
each containing exactly one index 1 and one index 2 critical point
of $f$. Then we can modify $\text{grad}(f)$ on these balls so that
we obtain a nonzero vector field on $M$ such that $v|\partial M =
v_0$. This shows that $\text{Spin}^c(M,\gamma) \neq \emptyset$. From
obstruction theory we get that $\text{Spin}^c(M, \gamma)$ is an
affine space over $H^2(M, \partial M; \mathbb{Z})$. 
\end{rem}

Next we define the Chern class of a $\text{Spin}^c$
structure. Let $i \co \partial M \to M$ denote the embedding. If
$v$ is the vector field constructed in \fullref{rem:3} then using
the naturality of Chern classes we see that $i^*(c_1(v^{\perp})) =
\delta$, thus $\delta \in \text{Im}(i^*)$. 


\begin{defn} \label{defn:12}
For $\mathfrak{s} \in \text{Spin}^c(M, \gamma)$ defined by a vector
field $v$ on $M$, let the class $c_1(\mathfrak{s}) = c_1(v^{\perp}) \in
(i^*)^{-1}(\delta)$ where $v^{\perp}$ is the oriented two-plane
field on $M$ perpendicular to $v$. Note that a priori we only know
that $c_1(v^{\perp}) \in H^2(M; \mathbb{Z})$, but since $v|\partial
M = v_0$ we get that $c_1(v^{\perp}) \in (i^*)^{-1}(\delta)$. 
\end{defn}

Fix a balanced diagram $(\Sigma, \boldsymbol{\alpha},
\boldsymbol{\beta})$ for $(M, \gamma)$. 

\begin{defn} \label{defn:13}
To each $\mathbf{x} \in \mathbb{T}_{\alpha} \cap \mathbb{T}_{\beta}$
we assign a $\text{Spin}^c$ structure $\mathfrak{s}(\mathbf{x}) \in
\text{Spin}^c(M, \gamma)$ as follows. Choose a Morse function $f$ on
$M$ compatible with the given balanced diagram $(\Sigma,
\boldsymbol{\alpha}, \boldsymbol{\beta})$. Then $\mathbf{x}$
corresponds to a multi-trajectory $\gamma_{\mathbf{x}}$ of
$\text{grad}(f)$ connecting the index one and two critical points of
$f$. In a regular neighborhood $N(\gamma_{\mathbf{x}})$ we can
modify $\text{grad}(f)$ to obtain a nonzero vector field $v$ on $M$
such that $v|\partial M = v_0$. We define $\mathfrak{s}(\mathbf{x})$
to be the homology class of this vector field $v$. 
\end{defn}

\begin{defn} \label{defn:14}
Let $\mathbf{x}, \mathbf{y} \in \mathbb{T}_{\alpha} \cap
\mathbb{T}_{\beta}$ and let $\gamma_{\mathbf{x}}$, 
$\gamma_{\mathbf{y}}$ be the corresponding multi-trajec\-tories,
thought of as one-chains in $M$. Then define $\epsilon(x,y) =
\gamma_{\mathbf{x}} - \gamma_{\mathbf{y}} \in H_1(M; \mathbb{Z})$. 

Alternatively, we can define $\epsilon(\mathbf{x}, \mathbf{y})$ in
the following manner. Choose paths $a \co I \to
\mathbb{T}_{\alpha}$ and $b \co I \to \mathbb{T}_{\beta}$ with
$\partial a = \partial b = \mathbf{x} - \mathbf{y}$. Then $a-b$ can
be viewed as a one-cycle in $\Sigma$ whose homology class in $M$ is
$\epsilon(\mathbf{x}, \mathbf{y})$. This is independent of the
choice of $a$ and $b$. 
\end{defn}

\begin{lem} \label{lem:5}
For $\mathbf{x}, \mathbf{y} \in \mathbb{T}_{\alpha} \cap
\mathbb{T}_{\beta}$ we have that $\mathfrak{s}(\mathbf{x})-
\mathfrak{s}(\mathbf{y}) = PD[\epsilon(\mathbf{x}, \mathbf{y})]$, 
where $PD \co H_1(M, \mathbb{Z}) \to H^2(M, \partial M;
\mathbb{Z})$ is the Poincar\'e duality map.
\end{lem}

\begin{proof}
The vector fields $\mathfrak{s}(\mathbf{x})$ and
$\mathfrak{s}(\mathbf{y})$ differ only in a neighborhood of
$\gamma_{\mathbf{x}}-\gamma_{\mathbf{y}}$. It is now a local
calculation to see that
$\mathfrak{s}(\mathbf{x})-\mathfrak{s}(\mathbf{y}) =
PD[\gamma_{\mathbf{x}}-\gamma_{\mathbf{y}}]$ (see \cite[Lemma
2.19]{OSz}).
\end{proof}

\begin{cor} \label{cor:4}
If $\mathfrak{s}(\mathbf{x}) \neq \mathfrak{s}(\mathbf{y})$ then
$D(\mathbf{x}, \mathbf{y}) = \emptyset$. In particular, there is no
Whitney disc connecting $\mathbf{x}$ and $\mathbf{y}$. If
$\pi_1(\text{Sym}^d(\Sigma)) = H_1(\text{Sym}^d(\Sigma);
\mathbb{Z})$ then the converse also holds, ie,
$\mathfrak{s}(\mathbf{x}) = \mathfrak{s}(\mathbf{y})$ implies that
$\pi_2(\mathbf{x},\mathbf{y}) \neq \emptyset$. 
\end{cor}

\begin{prop} \label{prop:6}
If $d > 1$ then $\pi_1(\text{Sym}^d(\Sigma)) =
H_1(\text{Sym}^d(\Sigma); \mathbb{Z})$. 
\end{prop}

\begin{proof}
The proof is analogous to \cite[Lemma 2.6]{OSz}. Let $\gamma\co
S^1 \to \text{Sym}^d(\Sigma)$ be a null-homologous curve that misses
the diagonal. Then there is a $2$--manifold with boundary $F$, a map
$i \co F \to \Sigma$ and a $d$--fold covering $\pi \co
\partial F \to S^1$ such that $(i|\partial F) \circ \pi^{-1} =
\gamma$. By increasing the genus of $F$ is necessary, we can extend
the $d$--fold covering of $S^1$ to a branched $d$--fold covering $\pi
\co F \to D^2$. Then the map $i \circ \pi^{-1} \co D^2 \to
\text{Sym}^d(\Sigma)$ shows that $\gamma$ is null-homotopic.
\end{proof}

\section{Maslov index}

Fix a balanced sutured manifold $(M,\gamma)$ and a balanced diagram
$(\Sigma, \boldsymbol{\alpha}, \boldsymbol{\beta})$ defining it.

\begin{note} \label{note:3}
For $\mathbf{x}, \mathbf{y} \in \mathbb{T}_{\alpha} \cap
\mathbb{T}_{\beta}$ and for a homotopy class $\phi \in
\pi_2(\mathbf{x}, \mathbf{y})$ let $\mathcal{M}(\phi)$ denote the
moduli space of pseudo-holomorphic representatives of $\phi$, and
let $\smash{\widehat{\mathcal{M}}}(\phi)$ be the quotient of this moduli
space by the action of $\mathbb{R}$. Let $\mu(\phi)$ denote the
Maslov index of $\phi$, ie, the expected dimension of
$\mathcal{M}(\phi)$. 
\end{note}

\begin{thm} \label{thm:1}
For $\mathbf{x} \in \mathbb{T}_{\alpha} \cap \mathbb{T}_{\beta}$ and
$\psi \in \pi_2(\mathbf{x},\mathbf{x})$ we have $$\mu(\psi) =
\langle c_1(\mathfrak{s}(\mathbf{x})), H(\mathcal{D}(\psi))
\rangle.$$
\end{thm}

\begin{proof}
See \cite[Theorem 4.9]{OSz}.
\end{proof}

\begin{cor} \label{cor:2}
Suppose that for $\phi_1, \phi_2 \in \pi_2(\mathbf{x},\mathbf{y})$
we have that $\mathcal{D}(\phi_1) = \mathcal{D}(\phi_2)$. Then
$\mu(\phi_1) = \mu(\phi_2)$. 
\end{cor}

\begin{proof}
The homotopy class $\psi = \phi_1 \overline{\phi}_2 \in
\pi_2(\mathbf{x},\mathbf{x})$ satisfies $\mathcal{D}(\psi) =
\mathcal{D}(\phi_1) - \mathcal{D}(\phi_2) = 0$ and $\mu(\psi) =
\mu(\phi_1) - \mu(\phi_2)$. The result then follows from \fullref{thm:1} using the fact that $H(\mathcal{D}(\psi)) = 0$. 
\end{proof}

This justifies the following definition.

\begin{defn} \label{defn:18}
We define the Maslov index of a domain $\mathcal{D} \in D(\Sigma,
\boldsymbol{\alpha}, \boldsymbol{\beta})$ as follows. If there exists a
homotopy class $\phi$ of Whitney discs such that $\mathcal{D}(\phi)
= \mathcal{D}$ then define $\mu(\mathcal{D}) := \mu(\phi)$. Otherwise we
define $\mu(\mathcal{D})$ to be $-\infty$. Furthermore, let
$\mathcal{M}(\mathcal{D})$ denote the moduli space of holomorphic
Whitney discs $u$ such that $\mathcal{D}(u) = \mathcal{D}$ and let
$\widehat{\mathcal{M}}(\mathcal{D}) =
\mathcal{M}(\mathcal{D})/\mathbb{R}$. 
\end{defn}

Thus we can rephrase \fullref{thm:1} as follows.

\begin{thm} \label{thm:2}
For $\mathbf{x} \in \mathbb{T}_{\alpha} \cap \mathbb{T}_{\beta}$ and
$\mathcal{P} \in D(\mathbf{x}, \mathbf{x})$ such that
$\mu(\mathcal{P}) \neq -\infty$ we have
$$\mu(\mathcal{P}) = \langle c_1(\mathfrak{s}(\mathbf{x})),
H(\mathcal{P}) \rangle.$$
\end{thm}

\section{Energy bounds}

First we recall the definition of the energy of a map of a planar
domain into a Riemannian manifold.

\begin{defn}
Let $\Omega$ be a domain in $\mathbb{C}$ and let $(X, g)$ be a
Riemannian manifold. The energy of a smooth map $u \co \Omega \to
X$ is given by $$E(u) = \frac12 \int_{\Omega} |du|^2.$$
\end{defn}

Let $\eta$ be a K\"ahler form on $\Sigma$. 

\begin{defn}
The \emph{area\/} of a domain $\mathcal{D} = \sum_{i=1}^m
n_i\mathcal{D}_i \in D(\Sigma, \boldsymbol{\alpha},
\boldsymbol{\beta})$ is defined as
$$\mathcal{A}(\mathcal{D}) = \sum_{i=1}^m n_i \cdot
\text{Area}_{\eta}(\mathcal{D}_i),
$$ where $\text{Area}_{\eta}(\mathcal{D}_i) = \int_{\mathcal{D}_i}\eta$. 
\end{defn}

\begin{thm} \label{thm:3}
There is a constant $C$ which depends only on the balanced diagram
$(\Sigma, \boldsymbol{\alpha}, \boldsymbol{\beta})$ and $\eta$ such
that for any smooth Whitney disc $$u \co (\mathbb{D},
\partial \mathbb{D}) \to (\text{Sym}^d(\Sigma), \mathbb{T}_{\alpha}
\cup \mathbb{T}_{\beta})$$ we have the energy bound
$$E(u) \le C \cdot \mathcal{A}(\mathcal{D}(u)).$$
\end{thm}

\begin{proof}
The proof is analogous to the proof of \cite[Lemma 3.5]{OSz}. We use
the fact that $\Sigma$ is compact. See also the paragraph below
\cite[Remark 3.7]{OSz}.
\end{proof}

\begin{cor} \label{cor:3}
For any $\mathcal{D} \in D(\Sigma, \boldsymbol{\alpha},
\boldsymbol{\beta})$ such that $\mu(\mathcal{D}) = 1$ the moduli
space $\widehat{\mathcal{M}}(\mathcal{D})$ is a compact
zero-dimensional manifold.
\end{cor}

\begin{proof}
This follows from \fullref{thm:3} using Gromov compactness.
\end{proof}

\begin{lem} \label{lem:8}
Every pseudo-holomorphic map $$u \co (\mathbb{D}, \partial
\mathbb{D}) \to (\text{Sym}^d(\Sigma), \mathbb{T}_{\alpha})$$ is
constant. The same holds for pseudo-holomorphic maps $$u' \co
(\mathbb{D}, \partial \mathbb{D}) \to (\text{Sym}^d(\Sigma),
\mathbb{T}_{\beta}).$$ Finally, every pseudo-holomorphic sphere $v
\co S^2 \to \text{Sym}^d(\Sigma)$ is constant.
\end{lem}

\begin{proof}
The boundary of the domain $\mathcal{D}(u)$ is a linear combination
of the $\alpha$ curves. Since $\alpha_1, \dots, \alpha_d$ are
linearly independent in $\Sigma$, this implies that $\mathcal{D}(u)
= 0$. Thus using \fullref{thm:3} we get that $E(u) \le C \cdot
\mathcal{A}(0) = 0$, so $E(u) =0$. A pseudo-holomorphic map with
zero energy is constant. A similar argument applies to $u'$. 

The domain of $v$ is a $2$--cycle. But $\Sigma$ has no closed
components, so $\mathcal{D}(v) = 0$. The fact that $v$ is constant
now follows similarly.
\end{proof}

\section{Definition of the chain complex}

Let $(M, \gamma)$ be a balanced sutured manifold and $(\Sigma,
\boldsymbol{\alpha}, \boldsymbol{\beta})$ an admissible balanced
diagram defining it. Fix a coherent system of orientations as in
\cite[Definition 3.11]{OSz}.

\begin{defn} \label{defn:19}
Let $CF(\Sigma, \boldsymbol{\alpha}, \boldsymbol{\beta})$ be the
free abelian group generated by the points in $\mathbb{T}_{\alpha}
\cap \mathbb{T}_{\beta}$. We define an endomorphism $\partial \co
CF(\Sigma, \boldsymbol{\alpha}, \boldsymbol{\beta}) \to CF(\Sigma,
\boldsymbol{\alpha}, \boldsymbol{\beta})$ so that for each generator
$\mathbf{x} \in \mathbb{T}_{\alpha} \cap \mathbb{T}_{\beta}$ we have
$$\partial \mathbf{x} = \sum_{\mathbf{y} \in
\mathbb{T}_{\alpha} \cap \,\mathbb{T}_{\beta}} \sum_{\{\mathcal{D}
\in D(\mathbf{x}, \mathbf{y})\,:\, \mu(\mathcal{D})=1\}}\#
\widehat{\mathcal{M}}(\mathcal{D})\cdot \mathbf{y}.
$$
\end{defn}

Since the diagram is admissible \fullref{lem:7} ensures that
$D(\mathbf{x}, \mathbf{y})$ has only finitely many positive
elements. But we know that from $\mathcal{M}(\mathcal{D}) \neq
\emptyset$ it follows that $\mathcal{D} \ge 0$. \fullref{cor:3} implies that if $\mu(\mathcal{D})=1$ then $\widehat{\mathcal{M}}(\mathcal{D})$ is a
compact zero-dimensional manifold, and the
coherent orientation system makes it oriented. Thus
$\#\widehat{\mathcal{M}}(\mathcal{D})$ makes sense, and the sum
above has only finitely many nonzero terms.

\begin{thm} \label{thm:4}
The pair $(CF(\Sigma, \boldsymbol{\alpha}, \boldsymbol{\beta}),
\partial)$ is a chain complex, ie, $\partial^2 =0$. 
\end{thm}

\begin{proof}
$\partial^2 = 0$ follows as in \cite[Theorem 4.1]{OSz}. Boundary
degenerations and spheres bubbling off are excluded by \fullref{lem:8}.
\end{proof}

\begin{defn}
For $\mathfrak{s} \in \text{Spin}^c(M, \gamma)$ let $C(\Sigma,
\boldsymbol{\alpha}, \boldsymbol{\beta}, \mathfrak{s})$ be the free
abelian group generated by those intersection points $\mathbf{x} \in
\mathbb{T}_{\alpha} \cap \mathbb{T}_{\beta}$ for which
$\mathfrak{s}(\mathbf{x}) = \mathfrak{s}$. 
\end{defn}

It follows from \fullref{cor:4} that $C(\Sigma,
\boldsymbol{\alpha}, \boldsymbol{\beta}, \mathfrak{s})$ is a
subcomplex of $C(\Sigma, \boldsymbol{\alpha}, \boldsymbol{\beta})$
and $$(C(\Sigma, \boldsymbol{\alpha}, \boldsymbol{\beta}),
\partial) = \bigoplus_{\mathfrak{s}\,
\in \text{Spin}^c(M, \gamma)} (C(\Sigma, \boldsymbol{\alpha},
\boldsymbol{\beta}, \mathfrak{s}),
\partial).$$
\begin{defn}
We define $\SFH(\Sigma, \boldsymbol{\alpha}, \boldsymbol{\beta})$ to
be the homology of the chain complex $(CF(\Sigma,
\boldsymbol{\alpha}, \boldsymbol{\beta}),
\partial)$. Similarly, for $\mathfrak{s} \in \text{Spin}^c(M,
\gamma)$ let $\SFH(\Sigma, \boldsymbol{\alpha}, \boldsymbol{\beta},
\mathfrak{s})$ be the homology of $(CF(\Sigma, \boldsymbol{\alpha},
\boldsymbol{\beta}, \mathfrak{s}),
\partial)$. 
\end{defn}

\begin{thm} \label{thm:5}
If the admissible balanced diagrams $(\Sigma, \boldsymbol{\alpha},
\boldsymbol{\beta})$ and $(\Sigma', \boldsymbol{\alpha}',
\boldsymbol{\beta}')$ define the same sutured manifold then
\begin{align*}\SFH(\Sigma, \boldsymbol{\alpha}, \boldsymbol{\beta}) &= \SFH(\Sigma',
\boldsymbol{\alpha}', \boldsymbol{\beta}')\notag\\
\SFH(\Sigma,
\boldsymbol{\alpha}, \boldsymbol{\beta}, \mathfrak{s}) &=
\SFH(\Sigma', \boldsymbol{\alpha}', \boldsymbol{\beta}',
\mathfrak{s}) \tag*{\hbox{and}}\end{align*} holds for every $\mathfrak{s} \in \text{Spin}^c(M,
\gamma)$. 
\end{thm}

\begin{proof}
This follows from \fullref{prop:4} as in \cite{OSz}.
\end{proof}

Thus we can make the following definition.

\begin{defn}
For $(M, \gamma)$ a balanced sutured manifold, we define the
\emph{sutured Floer homology\/} $\SFH(M, \gamma)$ as follows. Choose an
admissible balanced diagram $(\Sigma, \boldsymbol{\alpha},
\boldsymbol{\beta})$ defining $(M,\gamma)$. Then let $\SFH(M,\gamma)
= \SFH(\Sigma, \boldsymbol{\alpha, \boldsymbol{\beta}})$. For
$\mathfrak{s} \in \text{Spin}^c(M, \gamma)$ define $\SFH(M, \gamma,
\mathfrak{s})$ to be $\SFH(\Sigma, \boldsymbol{\alpha},
\boldsymbol{\beta}, \mathfrak{s})$. 
\end{defn}

\section{Relative gradings}

Suppose that $d>1$ in the balanced diagram $(\Sigma,
\boldsymbol{\alpha}, \boldsymbol{\beta})$. Then, using \fullref{prop:6} we get that $\pi_1(\text{Sym}^d(\Sigma)) =
H_1(\text{Sym}^d(\Sigma); \mathbb{Z})$. Thus, according to \fullref{cor:4}, for every $\mathbf{x}, \mathbf{y} \in
\mathbb{T}_{\alpha} \cap \mathbb{T}_{\beta}$ the equality
$\mathfrak{s}(\mathbf{x}) = \mathfrak{s}(\mathbf{y})$ implies that
$\pi_2(\mathbf{x},\mathbf{y}) \neq \emptyset$. Note that every
balanced sutured manifold $(M, \gamma)$ has a diagram with $d>1$, we
can achieve this by stabilizing an arbitrary balanced diagram
defining $(M, \gamma)$. 

\begin{defn}
For $\mathfrak{s} \in \text{Spin}^c(M,\gamma)$ let
$$\mathfrak{d}(\mathfrak{s}) = \gcd_{\xi \in H_2(M; \mathbb{Z})}
\langle c_1(\mathfrak{s}), \xi \rangle.
$$
\end{defn}

\begin{defn}
Let $\mathfrak{s} \in \text{Spin}^c(M,\gamma)$ and let $(\Sigma,
\boldsymbol{\alpha}, \boldsymbol{\beta})$ be an admissible balanced
diagram with $d>1$ defining $(M, \gamma)$. Then we define a relative
$\mathbb{Z}_{\mathfrak{d}(\mathfrak{s})}$ grading on $CF(\Sigma,
\boldsymbol{\alpha}, \boldsymbol{\beta}, \mathfrak{s})$ such that
for any $\mathbf{x}, \mathbf{y} \in \mathbb{T}_{\alpha} \cap
\mathbb{T}_{\beta}$ with $\mathfrak{s}(\mathbf{x}) =
\mathfrak{s}(\mathbf{y}) = \mathfrak{s}$ we have
$$\text{gr}(\mathbf{x}, \mathbf{y}) = \mu(\phi) \mod
\mathfrak{d}(\mathfrak{s}),
$$ where $\phi \in \pi_2(\mathbf{x}, \mathbf{y})$ is an arbitrary
homotopy class.
\end{defn}

The number $\text{gr}(\mathbf{x}, \mathbf{y})$ is independent of the
choice of $\phi$ because of \fullref{thm:2}. From the definition
of $\partial$ it is clear that $\text{gr}$ descends to a relative
grading on $\SFH(M, \gamma, \mathfrak{s})$. This grading is
independent of the balanced diagram defining the sutured manifold
$(M, \gamma)$. 

\section{Special cases and sample computations}

\begin{prop} \label{prop:7}
Let $Y$ be a closed connected oriented $3$--manifold. Then
$$\widehat{\HF}(Y) \approx \SFH(Y(1)).$$
\end{prop}

Recall that $Y(1)$ was introduced in \fullref{ex:1}. For the
definition of $\widehat{\HF}(Y)$ see \cite{OSz}.

\begin{proof}
Let $(\Sigma, \boldsymbol{\alpha}, \boldsymbol{\beta}, z)$ be a
weakly admissible Heegaard diagram defining $Y$. Choose a small
neighborhood $U$ of $z$ diffeomorphic to an open disc and let
$\Sigma' = \Sigma \setminus U$. Then $(\Sigma', \boldsymbol{\alpha},
\boldsymbol{\beta})$ is an admissible balanced sutured diagram
defining $Y(1)$. Since each $\mathcal{D} \in D(\Sigma',
\boldsymbol{\alpha}, \boldsymbol{\beta})$ has multiplicity zero at
$\partial \Sigma'$ the chain complexes $\widehat{CF}(\Sigma,
\boldsymbol{\alpha}, \boldsymbol{\beta},z)$ and $CF(\Sigma',
\boldsymbol{\alpha}, \boldsymbol{\beta})$ are isomorphic.
\end{proof}

In \fullref{ex:2} for every link $L$ in a closed connected
oriented $3$--manifold $Y$ we defined a balanced sutured manifold
$Y(L)$. In \cite{OSz2} an invariant $\smash{\widehat{\HFL}(\vec{L})}$ was
defined for oriented links $\vec{L} \subset Y$. Suppose that
$L$ has $l$ components, then $\smash{\widehat{\HFL}(\vec{L})}$ is computed
using $2l$--pointed Heegaard diagrams and Floer homology is taken
with coefficients in $\mathbb{Z}_2$. 

\begin{prop}
If $Y$ is a closed connected oriented $3$--manifold, $L \subset Y$ is a
link, and $\vec{L}$ is an arbitrary orientation of $L$ then
$$\widehat{\HFL}(\vec{L}) \approx \SFH(Y(L)) \otimes \mathbb{Z}_2.
$$ If $L$ has only one component, ie, if $L$ is a knot $K$, then
$$\widehat{\HFK}(Y, K) \approx \SFH(Y(K)).$$
\end{prop}

For the definition of $\widehat{\HFK}$ see \cite{OSz3} or \cite{Ras}.

\begin{proof}
Let $l$ be the number of components of the link $L$. If $(\Sigma,
\boldsymbol{\alpha}, \boldsymbol{\beta}, \mathbf{w}, \mathbf{z})$ is
a weakly admissible $2l$--pointed balanced Heegaard diagram of
$\vec{L}$ in the sense of \cite{OSz2} then remove an open regular
neighborhood of $\mathbf{w} \cup \mathbf{z}$ to obtain a compact
surface $\Sigma'$. The diagram $(\Sigma', \boldsymbol{\alpha},
\boldsymbol{\beta})$ is a balanced diagram defining the sutured
manifold $Y(L)$. It is now clear that the two chain complexes are
isomorphic.
\end{proof}

\begin{rem} \label{rem:4}
Suppose that $(M, \gamma)$ is a balanced sutured manifold such
that $\partial M$ is a torus and $s(\gamma)$ consists of two
components that represent the same homology class in $H_1(\partial
M; \mathbb{Z})$. Then $\SFH(M, \gamma)$ can be computed as the knot
Floer homology of the knot obtained from $(M, \gamma)$ using the
Dehn filling construction described in \fullref{ex:2}.
\end{rem}

\begin{prop} \label{prop:9}
If $(M, \gamma)$ is a product sutured manifold then
$$\SFH(M, \gamma) \approx \mathbb{Z}.$$
\end{prop}

\begin{proof}
Since $(M, \gamma)$ is product there is a compact oriented surface
$R$ with no closed components such that $(M, \gamma) = (R \times I,
\partial R \times I)$. Then $(R, \boldsymbol{\alpha},
\boldsymbol{\beta})$ is a balanced diagram defining $(M, \gamma)$,
where $\boldsymbol{\alpha} = \emptyset$ and $\boldsymbol{\beta} =
\emptyset$. Since $H_2(M; \mathbb{Z}) = 0$ any balanced diagram
defining $(M, \gamma)$ is admissible. Thus $\text{Sym}^0(\Sigma)
=\{\text{pt}\}$ and $\mathbb{T}_{\alpha} = \{\text{pt}\} =
\mathbb{T}_{\beta}$. Hence $\mathbb{T}_{\alpha} \cap
\mathbb{T}_{\beta}$ consists of a single point. Alternatively, we
can stabilize the above diagram and obtain the same result.
\end{proof}

\begin{rem}
Let $P$ denote the Poincar\'e $3$--sphere. Then the balanced sutured
manifold $P(1)$ is not a product. On the other hand $\SFH(P(1)) \approx \widehat{\HF}(P) \approx
\mathbb{Z}$ by \fullref{prop:7}. 
\end{rem}

\begin{defn}
A sutured manifold $(M, \gamma)$ is called \emph{irreducible\/} if
every $2$--sphere smoothly embedded in $M$ bounds a $3$--ball.
\end{defn}

\begin{qn} \label{qn:1}
Is the converse of \fullref{prop:9} true under certain
hypotheses? More precisely, suppose that the manifold $(M,\gamma)$ is irreducible
and $H_2(M; \mathbb{Z}) = 0$. Then does $\SFH(M, \gamma) =
\mathbb{Z}$ imply that $(M, \gamma)$ is a product sutured manifold?
\end{qn}

Next we recall the definition of a sutured manifold decomposition; see \cite[Definition 3.1]{Gabai}.

\begin{defn}
Let $(M, \gamma)$ be a sutured manifold and $S$ a properly embedded
oriented surface in $M$ such that for every component $\lambda$ of
$S \cap \gamma$, one of (1)-(3) holds:

(1)\qua $\lambda$ is a properly embedded nonseparating arc in $\gamma$. 

(2)\qua $\lambda$ is a simple closed curve in an annular component $A$
of $\gamma$ in the same homology class as $A \cap s(\gamma)$. 

(3)\qua $\lambda$ is a homotopically nontrivial curve in a torus
component $T$ of $\gamma$, and if $\delta$ is another component of
$T \cap S$, then $\lambda$ and $\delta$ represent the same homology
class in $H_1(T)$. 

Then $S$ defines a \emph{sutured manifold decomposition\/} $$(M,
\gamma)\rightsquigarrow^{S} (M', \gamma'),$$ where $M' = M \setminus
\text{Int}(N(S))$ and 
\begin{gather*}
\gamma' = (\gamma \cap M') \cup N(S'_+ \cap
R_-(\gamma)) \cup N(S'_- \cap R_+(\gamma)).\notag\\
R_+(\gamma') =
((R_+(\gamma) \cap M') \cup S'_+) \setminus \text{Int}(\gamma')\tag*{\hbox{Also,}}\\
R_-(\gamma') = ((R_-(\gamma) \cap M') \cup S'_-) \setminus
\text{Int}(\gamma'),\tag*{\hbox{and}}
\end{gather*}
where $S'_+$ ($S'_-$) is the component of
$\partial N(S) \cap M'$ whose normal vector points out of (into)
$M'$. 
\end{defn}


\begin{rem}
In other words the sutured manifold $(M', \gamma')$ is constructed
by splitting $M$ along $S$, creating $R_+(\gamma')$ by adding $S'_+$
to what is left of $R_+(\gamma)$ and creating $R_-(\gamma')$ by
adding $S'_-$ to what is left of $R_-(\gamma)$. Finally, one creates
the annuli of $\gamma'$ by ``thickening'' $R_+(\gamma') \cap
R_-(\gamma')$. 
\end{rem}

\begin{exm}
If $L \subset Y$ is a link and if $R$ is a Seifert surface of $L$
then there is a sutured manifold decomposition
$Y(L)\rightsquigarrow^R Y(R)$. Furthermore, if $L$ has $l$
components then there is a sutured manifold decomposition
$Y(l)\rightsquigarrow^A Y(L)$, where $A \subset Y(l)$ is a union of
embedded annuli ``around'' the link $L$. 
\end{exm}

The following definition can be found in \cite{Gabai2}.

\begin{defn}
A sutured manifold decomposition
$(M, \gamma) \rightsquigarrow^D (M', \gamma')$ where $D$ is a disc
properly embedded in $M$ and $|D \cap s(\gamma)| = 2$ 
is called a \emph{product decomposition\/}.
\end{defn}

\begin{rem}
If $(M, \gamma)$ is balanced and if $(M, \gamma) \rightsquigarrow^D
(M', \gamma')$ is a product decomposition then $(M',\gamma')$ is
also balanced.
\end{rem}

The following lemma will be very useful for computing sutured Floer
homology since we can simplify the topology of our sutured manifold
before computing the invariant.

\begin{lem} \label{lem:9}
Let $(M, \gamma)$ be a balanced sutured manifold. If $(M, \gamma)
\rightsquigarrow^D (M', \gamma')$ is a product decomposition then
$$\SFH(M, \gamma) = \SFH(M', \gamma').$$
\end{lem}

\begin{proof}
Let $N(D)$ be a regular neighborhood of $D$ and choose a
diffeomorphism $t \co N(D) \to [-1,4]^3$ mapping $D$ to
$\{3/2\} \times [-1,4]^2$ and sending $s(\gamma) \cap N(D)$ to $[-1,4]
\times \partial [-1,4] \times \{3/2\}$. Let $p_3 \co [-1,4]^3 \to
[-1,4]$ denote the projection onto the third factor. Then we can
extend the function $p_3 \circ t$ from $N(D)$ to a Morse function $f
\co M \to \mathbb{R}$ as described in the proof of \fullref{prop:2}. Note that $f$ has no critical points in $N(D)$ and
that $D$ is a union of flowlines of $\text{grad}(f)$ connecting
$R_-(\gamma)$ with $R_+(\gamma)$. From $f$ we obtain a balanced
diagram $(\Sigma, \boldsymbol{\alpha}, \boldsymbol{\beta})$ where
$\Sigma = f^{-1}(3/2)$. 

The arc $\delta = D \cap \Sigma$ has boundary on $\partial \Sigma$
and is disjoint from $\boldsymbol{\alpha}$ and $\boldsymbol{\beta}$. 
Since $\partial \delta \subset \partial \Sigma$ every domain
$\mathcal{D} \in D(\Sigma, \boldsymbol{\alpha}, \boldsymbol{\beta})$
has zero multiplicity in the domain containing $\delta$. Cutting
$\Sigma$ open along $\delta$ we obtain a surface $\Sigma'$. The
balanced diagram $(\Sigma', \boldsymbol{\alpha},
\boldsymbol{\beta})$ defines the sutured manifold $(M', \gamma')$. 
Using \fullref{prop:5} isotope $(\Sigma',
\boldsymbol{\alpha}, \boldsymbol{\beta})$ to obtain an admissible
diagram $(\Sigma', \boldsymbol{\alpha}, \boldsymbol{\beta}')$ of
$(M',\gamma')$. Then $(\Sigma, \boldsymbol{\alpha},
\boldsymbol{\beta}')$ is an admissible diagram of $(M, \gamma)$
since every periodic domain $\mathcal{P} \neq 0$ has zero
multiplicity in the domain containing $\delta$ and thus corresponds
to a periodic domain in $(\Sigma', \boldsymbol{\alpha},
\boldsymbol{\beta}')$, so it has both positive and negative
multiplicities. Thus we can suppose that both diagrams $(\Sigma,
\boldsymbol{\alpha}, \boldsymbol{\beta})$ and $(\Sigma',
\boldsymbol{\alpha}, \boldsymbol{\beta})$ are admissible.

Since every domain $\mathcal{D} \in D(\Sigma, \boldsymbol{\alpha},
\boldsymbol{\beta})$ has zero multiplicity in the domain containing
$\delta$, the chain complexes $CF(\Sigma, \boldsymbol{\alpha},
\boldsymbol{\beta})$ and $CF(\Sigma', \boldsymbol{\alpha},
\boldsymbol{\beta})$ are isomorphic.
\end{proof}

As an application we prove a generalization of \fullref{prop:7}.

\begin{prop} \label{prop:13}
If $Y$ is a closed connected oriented $3$--manifold then for all $n
\ge 1$, $$\SFH(Y(n)) \approx \bigoplus_{2^{n-1}} \widehat{\HF}(Y)
\approx \widehat{\HF}(Y) \otimes \bigotimes_{n-1} \mathbb{Z}^2 .$$
\end{prop}

\begin{proof}
We prove the claim by induction on $n$. The case $n=1$ is true
according to \fullref{prop:7}. Suppose that we know
the proposition for some $n-1 \ge 1$. Then applying the induction
hypotheses to $(Y\#(S^1 \times S^2))(n-1)$ we get that $$\SFH
\left((Y\#(S^1 \times S^2))(n-1)\right) \approx
\widehat{\HF}\left(Y\#(S^1 \times S^2)\right) \otimes
\bigotimes_{n-2} \mathbb{Z}^2 \approx\widehat{\HF}(Y)
\otimes \bigotimes_{n-1} \mathbb{Z}^2.$$ Here we used the connected
sum formula $\widehat{\HF}(Y\#(S^1 \times S^2)) \approx
\widehat{\HF}(Y) \otimes \widehat{\HF}(S^1 \times S^2)$ and the fact
that $\widehat{\HF}(S^1 \times S^2) \approx \mathbb{Z}^2$. On the
other hand we will show that there is a product
decomposition $Y(n-1)\#(S^1 \times S^2) \rightsquigarrow^D Y(n)$, 
which shows together with \fullref{lem:9} that the induction
hypothesis is also true for $n$. 

To find the product disc $D$ choose a ball $B_1 \subset S^1 \times
S^2$ such that there is a point $p \in S^1$ for which $\{p\} \times
S^2$ intersects $B_1$ in a disc. Then let $D$ be the closure of
$(\{p\} \times S^2) \setminus B_1$. We can also choose a simple
closed curve $s_1 \subset \partial B_1$ so that $|s_1 \cap D| = 2$. 
Now construct $(Y\#(S^1 \times S^2))(n-1)$ as in \fullref{ex:1} using $B_1$ and $s_1$ as above . Then $D$ is a product
disc with the required properties.
\end{proof}

Next we will generalize the above idea to obtain a connected sum
formula.

\begin{prop} \label{prop:12}
Let $(M, \gamma)$ and $(N, \nu)$ be balanced sutured manifolds and
let $Y$ be a closed oriented $3$--manifold. Then \begin{align*}\SFH((M, \gamma)\#
(N, \nu)) &= \SFH(M, \gamma) \otimes \SFH(M, \nu) \otimes \mathbb{Z}^2.\notag\\
\SFH(M\# Y, \gamma) &= \SFH(M,\gamma) \otimes \widehat{\HF}(Y). 
\tag*{\hbox{Furthermore,}}\end{align*}
\end{prop}

\begin{proof}
There are product decompositions \begin{align*}(M, \gamma)\#(N, \nu)
&\rightsquigarrow^D (M, \gamma) \coprod N(1)\notag\\ (M, \gamma)\#Y
&\rightsquigarrow^D (M, \gamma) \coprod Y(1). \tag*{\hbox{and}}\end{align*} To see this push some
part of the boundary of $M$ containing a segment of $\gamma$ into
the connected sum tube using a finger move and repeat the idea
described in the proof of \fullref{prop:13} (also see \fullref{fig:2}). 

\begin{figure}[ht!]
\labellist
\small\hair 2pt
\pinlabel $M$ at 189 428 
\pinlabel $s$ [tr] at 243 433
\pinlabel $s$ [tl] at 261 417
\pinlabel $D$ at 333 427
\pinlabel {$Y \hbox{ or } N$} at 413 431
\endlabellist
\centering
\includegraphics{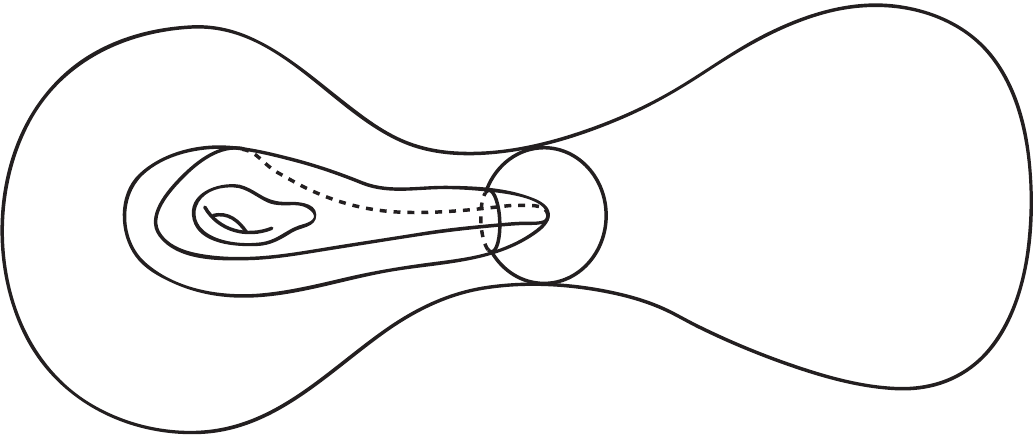}
\caption{Product decomposition of a connected sum} \label{fig:2}
\end{figure}
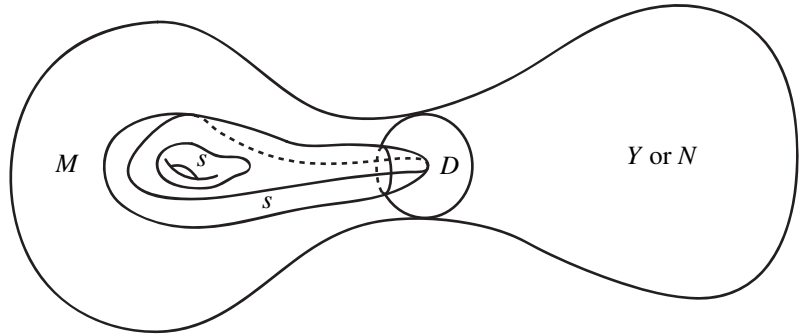

This implies that
\begin{gather} \label{eqn:1}
\SFH((M, \gamma)\#(N,\nu)) = \SFH(M, \gamma) \otimes \SFH(N(1))\\
\SFH(M\#Y, \gamma) = \SFH(M, \gamma) \otimes \SFH(Y(1)).\tag*{\hbox{and}}
\end{gather}
\fullref{prop:7} says that $\SFH(Y(1)) = \widehat{\HF}(Y)$. 
Since $N(1) = (N,\nu)\#S^3(1)$ we can apply \eqref{eqn:1}
again and we get that $$\SFH(N(1)) = \SFH(N, \nu) \otimes
\SFH(S^3(2)).$$ From the existence of a product decomposition $(S^1
\times S^2)(1) \rightsquigarrow S^3(2)$ (see \fullref{prop:13}) we obtain that
$$\SFH(S^3(2)) \approx \widehat{\HF}(S^1 \times S^2) \approx
\mathbb{Z}^2.\proved$$
\end{proof}

\begin{cor}
If $(M, \gamma)$ is a connected balanced sutured manifold and $n \ge
1$ then
$$\SFH(M(n)) \approx \SFH(M, \gamma) \otimes \mathbb{Z}^{2^n}.$$
\end{cor}

\begin{proof}
The claim follows by induction on $n$. The case $n=0$ is trivial.
Now let us suppose that $n>0$. Since $M(n) = M(n-1)\#S^3(1)$ we
get from \fullref{prop:12} that $\SFH(M(n)) \approx
\SFH(M(n-1)) \otimes \SFH(S^3(1)) \otimes \mathbb{Z}^2$. Here
$\SFH(S^3(1)) \approx \widehat{\HF}(S^3) \approx \mathbb{Z}$. This
concludes the proof.
\end{proof}

\begin{defn}
A sutured manifold $(M, \gamma)$ is called \emph{taut\/} if $M$ is
irreducible and $R(\gamma)$ is incompressible and Thurston
norm-minimizing in its homology class in  $H_2(M, \gamma)$. 
\end{defn}

\begin{prop} \label{prop:8}
Suppose that $(M, \gamma)$ is an irreducible balanced sutured
manifold. If $(M, \gamma)$ is not taut then $\SFH(M, \gamma) = 0$. 
\end{prop}

The following proof is due to Yi Ni.

\begin{proof}
Since $(M, \gamma)$ is not taut and $M$ is irreducible either
$R_+(\gamma)$ or $R_-(\gamma)$, say $R_+(\gamma)$, is either
compressible or it is not Thurston norm minimizing in $H_2(M,
\gamma)$. In both cases there exists a properly embedded surface
$(S,\partial S) \subset (M, \gamma)$ such that $\chi(S) >
\chi(R_+(\gamma))$, no collection of components of $S$ is
null-homologous and the class $[S,\partial S] = [R_+(\gamma), \partial
R_+(\gamma)]$ in $H_2(M, \gamma)$. Then decomposing $(M, \gamma)$
along $S$ we get two connected sutured manifolds $(M_+, \gamma_+)$
and $(M_-,\gamma_-)$. Here $R_+(\gamma) = R_+(\gamma_+)$ and
$R_-(\gamma) = R_-(\gamma_-)$. 

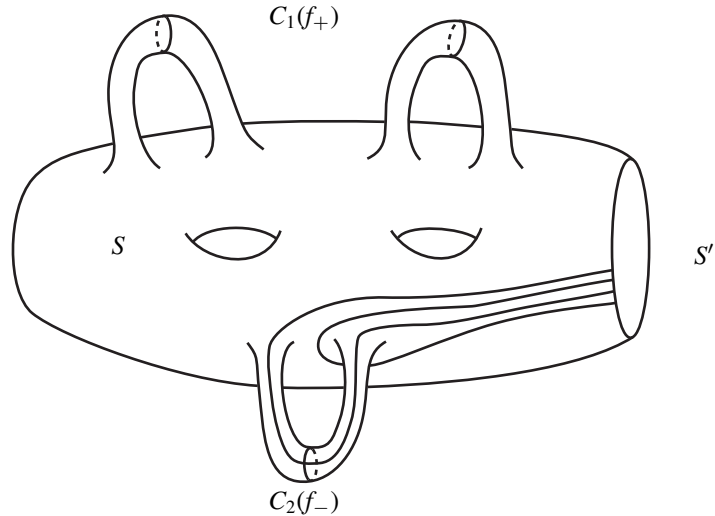
\begin{figure}[ht!]
\labellist
\small\hair 2pt
\pinlabel $S$ at 226 493
\pinlabel $S'$ at 447 490
\pinlabel {$C_{1}(f_{+})$} at 296 579
\pinlabel {$C_{2}(f_{-})$} at 296 396
\endlabellist
\centering
\includegraphics{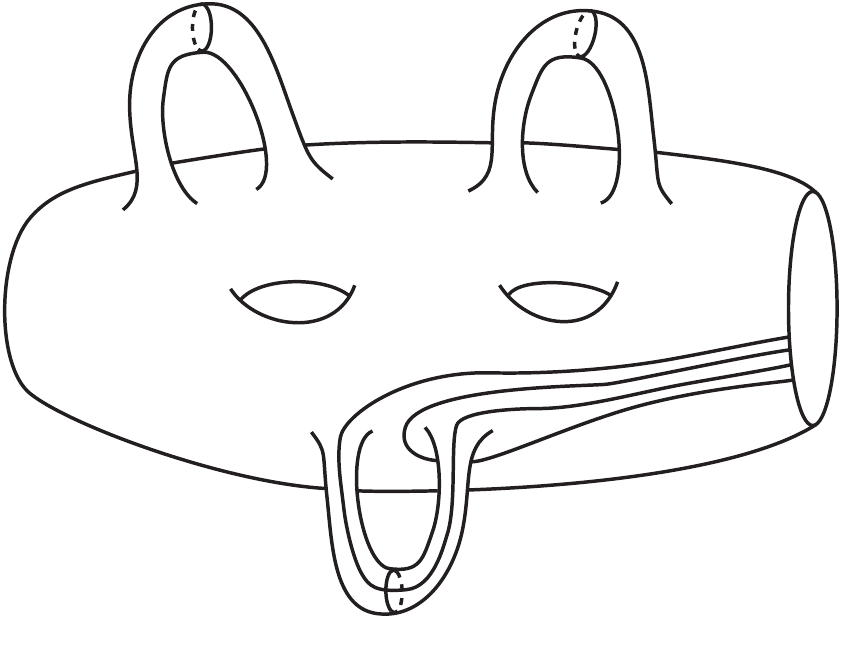}
\caption{The surface $S'$ with two $\boldsymbol{\alpha}_+$ curves on
the top, one $\boldsymbol{\beta_-}$ curve on the bottom, and with
two winding arcs, one of them intersecting $\boldsymbol{\beta}_-$
and being disjoint from $\boldsymbol{\alpha}_+$} \label{fig:3}
\end{figure}

As in the proof of \fullref{prop:2} construct Morse
functions $f_+$ and $f_-$ on $M_+$ and $M_-$, respectively, having no
index zero and three critical points. Then $f = f_+ \cup f_-$ is a
Morse function on $M$ that has $S$ as a level surface. Denote by $C_i(h)$ the set  
and by $c_i(h)$ the number of
index $i$ critical points of a Morse function $h$. Now rearrange $f$
by switching $C_1(f_+)$ and $C_2(f_-)$ to obtain a self-indexing
Morse function $g$ (see Milnor \cite{Milnor}). Then $g$ induces a Heegaard
diagram $(S', \boldsymbol{\alpha}_+ \cup \boldsymbol{\alpha}_-,
\boldsymbol{\beta}_+ \cup \boldsymbol{\beta}_-)$, where
$\boldsymbol{\alpha}_{\pm}$ and $\boldsymbol{\beta}_{\pm}$ are the
sets of attaching circles corresponding to the critical points in
$C_1(f_{\pm})$ and $C_2(f_{\pm})$, respectively, and $S'$ is
obtained by performing $c_1(f_+) + c_2(f_-)$ zero surgeries on $S$
whose belt circles are the elements of $\boldsymbol{\alpha}_+ \cup
\boldsymbol{\beta}_-$ (see \fullref{fig:3}). Our main observation
is that $\alpha \cap \beta = \emptyset$ if $\alpha \in
\boldsymbol{\alpha_+}$ and $\beta \in \boldsymbol{\beta_-}$, because
they are belt circles of two disjoint handles added to $S$. This
property of the Heegaard diagram is preserved if we apply the
winding argument of \fullref{prop:5} using winding arcs
$\gamma_1, \dots, \gamma_l$ that satisfy the following property: if
$\gamma_k \cap \beta \neq \emptyset$ for $\beta \in
\boldsymbol{\beta}_-$ and $1 \le k \le l$ then $\gamma_k \cap \alpha
= \emptyset$ for every $\alpha \in \boldsymbol{\alpha}_+$. Such arcs
$\gamma_1, \dots, \gamma_l$ are easy to construct (see \fullref{fig:3}). Thus we can assume our Heegaard diagram is admissible.

The only $\beta$ curves that can intersect $\alpha \in
\boldsymbol{\alpha}_+$ are the elements of $\boldsymbol{\beta}_+$. 
But $$0 > \chi(R_+(\gamma)) - \chi(S) = 2(c_2(f_+)-c_1(f_+)),$$ thus
$|\boldsymbol{\alpha}_+| = c_1(f_+) > c_2(f_+) =
|\boldsymbol{\beta}_+|$. This shows that $\mathbb{T}_{\alpha} \cap
\mathbb{T}_{\beta} = \emptyset$ for this Heegaard diagram. Indeed,
if there was a permutation $\pi \in S_d$ such that $\alpha_i \cap
\beta_{\pi(i)} \neq \emptyset$ for every $1 \le i \le d$, then for
$\alpha_i \in \boldsymbol{\alpha}_+$ we would have $\beta_{\pi(i)}
\in \boldsymbol{\beta}_+$, and the injectivity of $\pi$ would imply
that $|\boldsymbol{\beta}_+| \ge |\boldsymbol{\alpha}_+|$. 
\end{proof}

\begin{qn} \label{qn:2}
Is the converse of \fullref{prop:8} true, ie, if $\SFH(M, \gamma) = 0$
does it follow that $(M,\gamma)$ is not taut?
\end{qn}

\section{Seifert surfaces}

Now we turn our attention to \fullref{ex:4}. These sutured
manifolds are of particular interest to us due to the following
theorem of Gabai \cite[Theorem 1.9]{Gabai2}.

\begin{thm} \label{thm:6}
Suppose that $R$ is an oriented surface in $S^3$ and let $L$ be the
oriented link $\partial R$. Then $L$ is a fibred link with fibre $R$
if and only if $S^3(R)$ is a product sutured manifold.
\end{thm}

This becomes interesting in light of the following conjecture.

\begin{conj} \label{conj:1}
Let $K$ be a knot in $S^3$ and let $R$ be a genus $g$ Seifert
surface of $K$. Then $\widehat{\HFK}(K, g) \approx \SFH(S^3(R))$. 
\end{conj}

Note that from Alexander duality we get that $$H_2 \left(S^3(R);
\mathbb{Z}\right) \approx \widetilde{H}^0(R \times I; \mathbb{Z}) =
0.$$ Thus, together with a positive answer to \fullref{qn:2},
\fullref{conj:1} would give a new proof of the fact that
$\widehat{\HFK}(K,g(K)) \neq 0$, where $g(K)$ denotes the three-genus
of $K$. Combining \fullref{conj:1} with \fullref{prop:8} we would get that $\widehat{\HFK}(K,g) = 0$ for $g >
g(K)$. 

Finally, if we combine \fullref{thm:6}, \fullref{conj:1}
and \fullref{qn:1} we would obtain a proof of the following
conjecture (see \cite[Theorem 1.1]{OSz4} and \cite{Yi}). Note that a
fibred knot has a unique minimal genus Seifert surface up to
isotopy.

\begin{conj} \label{conj:2}
Let $K$ be a knot in $S^3$. Then $K$ is fibred if and only if
$$\widehat{\HFK}(K, g(K)) \approx \mathbb{Z}.$$
\end{conj}

In what follows we collect some evidence supporting \fullref{conj:1}. First we recall a result of Hedden \cite{Hedden}.

\begin{prop} \label{prop:10}
Let $K$ be knot in $S^3$ and let $D_+(K,t)$ denote the positive
$t$--twisted Whitehead double of $K$. The meridian $\mu \subset S^3$
of $K$ can be viewed as a knot in $S^3_t(K)$ (the parameter $t$ Dehn
surgery on $K$). Then we have $$\widehat{\HFK}(D_+(K,t),1) \approx
\widehat{\HFK}(S^3_t(K), \mu).
$$
\end{prop}

Using this result we can prove the following.

\begin{thm} \label{thm:7}
Let $R$ be the Seifert surface of $D_+(K,t)$ obtained by taking the
satellite of the surface $R_{2t}$ defined in \fullref{prop:11} (see below). Then
$$\widehat{\HFK}(D_+(K,t),1) \approx \SFH(S^3(R)).$$
\end{thm}

Note that the Seifert genus of $D_+(K,t)$ is $1$. The left hand side
of \fullref{fig:2} shows $D_+(U,1)$ together with its natural
Seifert surface. The surface $R$ is obtained by taking a solid torus
neighborhood of $U$ containing $R_{2t}$ and wrapping it around $K$
using the Seifert framing of $K$.

\begin{figure}[ht!]
\labellist
\small\hair 2pt
\pinlabel $\Tw_2$ at 145 398
\pinlabel $\Tw_3$ at 321 398
\pinlabel $D$ at 297 353
\endlabellist
\centering
\includegraphics[scale=.90]{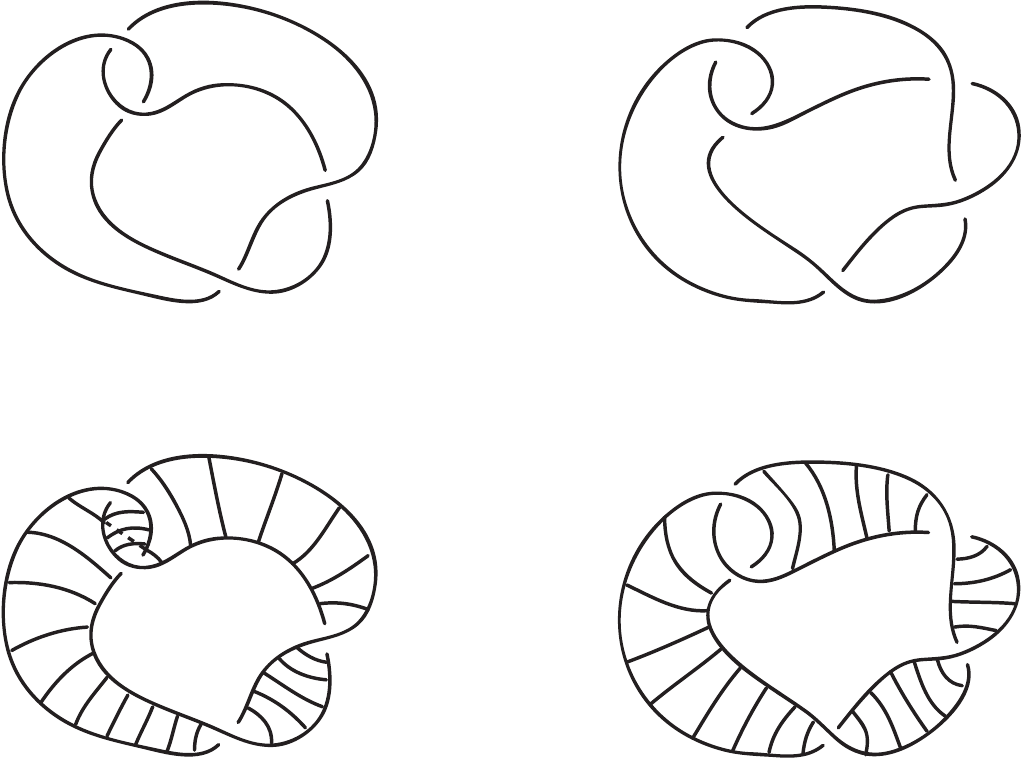}
\caption{Standard diagrams of the knots $\Tw_2 = D_+(U,1)$ and
$\Tw_3$ together with the Seifert surfaces $R_2$ and $R_3$
obtained using the Seifert\break algorithm} \label{fig:4}
\end{figure}

\begin{proof}
In light of \fullref{prop:10} we only have to show that
$$\SFH(S^3(R)) \approx \widehat{\HFK}(S^3_t(K), \mu).$$ Let $K_{2,2t}$
denote the $(2,2t)$--cable of $K$ (which is a two component link) and
let $R'$ be the natural Seifert surface of $K_{2,2t}$. Then there is
a product decomposition $S^3(R)\rightsquigarrow^D S^3(R')$ (see
\fullref{fig:4}). This does not change the sutured Floer homology
according to \fullref{lem:9}. Now we can apply \fullref{rem:4}
to compute $\SFH(S^3(R'))$. If we glue $S^1 \times D^2$ to $S^3(R')$
the meridian $\{1\} \times \partial D^2$ maps to one component of
$K_{2,2t}$ and we can suppose that the longitude $S^1 \times \{1\}$
maps to the meridian $\mu$ of the original knot $K$. After gluing in
$S^1 \times D^2$ we obtain $S^3_t(K)$. Note that in $S^3_t(K)$ the
knot $S^1 \times \{0\}$ is isotopic to $\mu$ since the longitude of
$S^1 \times D^2$ was identified with $\mu$. Thus $$\SFH(S^3(R'))
\approx \widehat{\HFK}(S^3_t(K), \mu),$$ which concludes the proof.
\end{proof}

The following similar statement can be proved without making use of
\fullref{prop:10}.

\begin{prop} \label{prop:11}
Let $\Tw_n$ denote the standard diagram of the twist knot with
$n$ half right-handed twists and let $R_n$ be the genus one Seifert
surface of $\Tw_n$ obtained from the Seifert algorithm (see
\fullref{fig:4}). Then
$$ \SFH(S^3(R_n)) \approx \widehat{\HFK}(\Tw_n,1) \approx
\mathbb{Z}^{[(n+1)/2]}.
$$
\end{prop}

\begin{proof}
Let $R_n'$ be the unique Seifert surface of the torus link
$T_{2,2[(n+1)/2]}$. Then there is a product decomposition
$S^3(R_n)\rightsquigarrow^D S^3(R_n')$. As in the proof of \fullref{thm:7} we have an isomorphism $$\SFH(S^3(R_n')) \approx
\widehat{\HFK}\left(S^3_{[(n+1)/2]}(U), \mu \right).$$ But
$S^3_{[(n+1)/2]}(U)$ is homeomorphic to the lens space
$L([(n+1)/2],1)$. Thus, according to \cite{OSz5},
$$\SFH(S^3(R_n)) \approx \widehat{\HFK}\left(L([(n+1)/2]), \mu \right)
\approx \mathbb{Z}^{[(n+1)/2]}.$$
Since $\Tw_n$ is alternating,
$\text{rk}\big(\widehat{\HFK}(\Tw_n,1)\big)$ agrees with the
absolute value of the leading coefficient of the Alexander
polynomial of $K$, which is $[(n+1)/2]$. 
\end{proof}

Finally one more evidence supporting \fullref{conj:1}.

\begin{prop}
Suppose that the knot $K$ has at most $7$ crossings and that $K \neq
7_4$. Then $K$ has a unique minimal genus Seifert surface $R$ and
$$\widehat{\HFK}(K, g(K)) \approx \SFH(S^3(R)).$$
\end{prop}

\begin{proof}
The fact that $K$ has a unique minimal genus Seifert surface was
proved by Kobayashi \cite{Kob}. We already know the statement for fibred
knots. The only nonfibred at most 7 crossing knots are $5_2, 6_1,
7_2, 7_3, 7_4$ and $7_5$. The knots $5_2, 6_1$ and $7_2$ are twist
knots and hence the result follows from \fullref{prop:11}.
The case of $7_3$ and $7_5$ is analogous, we can reduce the
computation of $\SFH(S^3(R))$ using product decompositions to
computing knot Floer homology of knots in lens spaces. Both knots
are alternating, so their knot Floer homology can be computed from
the Alexander polynomial.
\end{proof}

\begin{rem}
By understanding how a balanced diagram changes under a disc
decomposition of the underlying sutured manifold we could prove the
following formula. If the oriented surface $R \subset S^3$ is the
Murasugi sum of the surfaces $R_1$ and $R_2$ then over any field
$\mathbb{F}$ $$\SFH(S^3(R); \mathbb{F}) \approx \SFH(S^3(R_1);
\mathbb{F}) \otimes \SFH(S^3(R_2); \mathbb{F}).$$ This formula is
analogous to the Murasugi sum formula of \cite{Yi}.

The knot $7_4$ has two distinct minimal genus Seifert surfaces, both
of them Murasugi sums of two embedded annuli. In both cases we get
that the sutured Floer homology associated to the Seifert surface is
isomorphic to $\mathbb{Z}^4$, supporting \fullref{conj:1}. I
will deal with these results in a separate paper.
\end{rem}

\bibliographystyle{gtart}
\bibliography{link}

\end{document}